\documentclass[oneside, reqno] {amsart}

\usepackage{xypic}
\input xy
\xyoption{all}
\usepackage{epsfig}
\usepackage{amsthm}
\usepackage{amssymb}
\usepackage{amsmath}
\usepackage{amscd}
\usepackage{color}
\usepackage[T1]{fontenc}
\usepackage{upgreek}
\usepackage[font=scriptsize]{caption}
\usepackage{wrapfig}
\usepackage{graphicx}
\usepackage{subfigure}

\newcommand{\bg}{\begin{equation}}
\newcommand{\ed}{\end{equation}}
\newcommand{\bga}{\begin{eqnarray}}
\newcommand{\eda}{\end{eqnarray}}
\newcommand{\pf}{\textbf{Proof:\ }}

\def\cbdu{\par{\raggedleft$\Box$\par}}

\newtheorem {Theorem}  {Theorem}

\numberwithin{Theorem}{section}

\newtheorem {Lemma}[Theorem]  {Lemma}
\newtheorem {Proposition}[Theorem]{Proposition}
\theoremstyle{definition}
\newtheorem{Definition}[Theorem]{Definition}
\theoremstyle{remark}
\newtheorem{Remark}[Theorem]{\bf Remark}

\expandafter\chardef\csname pre amssym.def
at\endcsname=\the\catcode`\@ \catcode`\@=11
\def\undefine#1{\let#1\undefined}
\def\newsymbol#1#2#3#4#5{\let\next@\relax
 \ifnum#2=\@ne\let\next@\msafam@\else
 \ifnum#2=\tw@\let\next@\msbfam@\fi\fi
 \mathchardef#1="#3\next@#4#5}
\def\mathhexbox@#1#2#3{\relax
 \ifmmode\mathpalette{}{\m@th\mathchar"#1#2#3}%
 \else\leavevmode\hbox{$\m@th\mathchar"#1#2#3$}\fi}
\def\hexnumber@#1{\ifcase#1 0\or 1\or 2\or 3\or 4\or 5\or 6\or 7\or 8\or
 9\or A\or B\or C\or D\or E\or F\fi}

\font\teneufm=eufm10 \font\seveneufm=eufm7 \font\fiveeufm=eufm5
\newfam\eufmfam
\textfont\eufmfam=\teneufm \scriptfont\eufmfam=\seveneufm
\scriptscriptfont\eufmfam=\fiveeufm

\catcode`\@=\csname pre amssym.def at\endcsname

\newcounter{remark}
\newcommand{\supp}{{\mathit supp}\,}

\renewcommand{\div}{\mbox{div}}

\def  \12  {{\frac{1}{2}}}

\def\build#1_#2^#3{\mathrel{\mathop{\kern 0pt#1}\limits_{#2}^{#3}}}

\numberwithin{equation}{section}

\begin{document}
%\currannalsline{0}{2006}

\title[Forced SQG]{Non-unique weak solutions of forced SQG}

%\author{hello}

\author [Mimi Dai]{Mimi Dai}

\address{Department of Mathematics, Statistics and Computer Science, University of Illinois at Chicago, Chicago, IL 60607, USA}
\email{mdai@uic.edu}

\author [Qirui Peng]{Qirui Peng}

\address{Department of Mathematics, Statistics and Computer Science, University of Illinois at Chicago, Chicago, IL 60607, USA}
\email{qpeng9@uic.edu}

\thanks{The authors are partially supported by the NSF grants DMS--2009422 and DMS--2308208. } 
%M. Dai is also supported by the AMS Centennial Fellowship.}

\begin{abstract}

We construct non-unique weak solutions $\theta\in C_t^0C_x^{0-}$ for forced surface quasi-geostrophic (SQG) equation. This is achieved through a convex integration scheme adapted to the sum-difference system of two distinct solutions. Without external forcing, non-unique weak solutions $\theta$ in space $C_t^0C_x^{\alpha}$ with $\alpha<-\frac15$ were constructed by Buckmaster, Shkoller and Vicol \cite{BSV} and Isett and Ma \cite{IM}.

\bigskip


\hspace{0.02cm}CLASSIFICATION CODE: 35Q35, 35Q86, 76D03.
\end{abstract}

\maketitle

\section{Introduction}
\label{sec-int}

\subsection{Overview}

We consider the surface quasi-geostrophic equation (SQG) with external forcing
\begin{equation}\label{sqg}
\begin{split}
\partial_t\theta+ u \cdot\nabla \theta+\Lambda^\gamma \theta=&\ f, \\
u=&\ \nabla^{\perp}\Lambda^{-1}\theta %=&\left( -\mathcal R_2\theta, \mathcal R_1\theta\right)
\end{split}
\end{equation}
on $\mathbb T^2\times [0,\infty)$. 
In (\ref{sqg}), $\theta$ denotes the surface temperature in a rapidly rotating and stratified flow;  $u$ is velocity field of the flow; and $f$ represents the buoyancy forcing. The parameter $\gamma\geq 0$ indicates the strength of dissipation. By convention, $\Lambda^0 \theta=0$ when $\gamma=0$. The operator $\nabla^{\perp}\Lambda^{-1}$ with $\Lambda=(-\Delta)^{\frac12}$ is the Riesz transform, which is a non-local odd operator.
%The Zygmund operator $\Lambda$ is defined as $\Lambda=(-\Delta)^{\frac12}$; $R_1$ and $R_2$ are Riesz transforms. We assume $0<\gamma<\frac32$. It belongs to the family of active scalar equations with the non-local operator $T=: \nabla^{\perp}\Lambda^{-1}$ and drift velocity $u=T[\theta]$. Note the operator $T$ is odd in the sense that its Fourier symbol is odd; and $\nabla\cdot u= 0$.

As an important model in atmosphere and oceanography investigations and due to its analogous features with the 3D Euler equation, SQG has been studied to a great extent through the lens of rigorous mathematical analysis, tracing back to \cite{CMT, Ped}. Many survey papers on the development of SQG can be found in the literature, for instance, see the recent work \cite{BSV} and references therein. 

One significant property of the inviscid SQG without external forcing is the infinitely many conservation laws, including the Hamiltonian $\int |\Lambda^{-\frac12}\theta|^2\,dx$ and $L^p$ norms for $1\leq p\leq \infty$. Another crucial feature relies on the special structure of the nonlinearity $u\cdot\nabla\theta$ thanks to the Riesz transform, which leads to certain cancelation in analysis. These properties are beneficial in the study of regularity, well-posedness and long time dynamics problems for SQG. 
In particular, the author of \cite{Res} exploited the aforementioned cancelation property and showed 
global existence of weak solutions in $L^2$ to the unforced SQG with $0\leq\gamma\leq 2$, with the help of weak compactness. 
Moreover, existence of global regular solution for the critical SQG (\ref{sqg}) with $\gamma=1$ and $f\equiv 0$ was established by \cite{CaV, CV, KNV} independently via different sophisticated methods, which all took advantage of the special nonlinear structure. For the critical SQG (\ref{sqg}) with $\gamma=1$ in the presence of external forcing $f$, the authors of \cite{CTV1} showed the absence of anomalous dissipation of energy provided $f\in L^1(\mathbb R^2)\cap L^\infty(\mathbb R^2)$. Also for the forced critical SQG, global attractor was obtained in \cite{CTV2} for $f\in L^\infty(\mathbb T^2)\cap H^1(\mathbb T^2)$ and in \cite{CD} while $f\in L^p(\mathbb T^2)$ with $p>0$.

Recent attention has been centered around the generalized Onsager conjecture and unique solvability for SQG following the flourish of convex integration techniques in PDEs (c.f. \cite{BV, DLS1, DLS2}) and the resolution of Onsager's conjecture for the Euler equation \cite{BDLSV, Is}. For the inviscid SQG without external forcing, the authors of \cite{BSV} obtained non-unique weak solutions $\theta$ with $\Lambda^{-1}\theta\in C_t^\sigma C_x^\alpha$ for $\frac12<\alpha<\frac 45$ and $\sigma<\frac{\alpha}{2-\alpha}$. Note the Hamiltonian is expected to be conserved for solutions with regularity $\theta\in C^0_{t,x}$ and the constructed solutions in \cite{BSV} dissipates the Hamiltonian. By exploiting a nifty bilinear form associated with the nonlinear structure of SQG in \cite{IM}, the authors provided another elegant proof of non-uniqueness of weak solutions in the same regularity space. For the stationary SQG with dissipation $\Lambda^\gamma \theta$, the work \cite{CKL} gave a construction of non-trivial weak solutions $\theta$ satisfying $\Lambda^{-1}\theta\in C^\alpha$ for $\frac12<\alpha<\frac 23$ and $0<\gamma< 2-\alpha$. While for the forced stationary SQG, we \cite{DP} showed the non-uniqueness with  $\Lambda^{-1}\theta\in C^\alpha$ for $\frac12<\alpha<\frac 34$ and $0<\gamma< 2-\alpha$. The later result indicates the relaxation of the problem by allowing the presence of external forcing. Such relaxation is also seen for the forced Euler equation in \cite{BHP} where non-unique weak solutions with regularity higher than Onsager's $\frac13$ regularity were constructed.

The purpose of this paper is to construct non-unique weak solutions for the forced SQG (\ref{sqg}) with regularity approaching the expected threshold $C^0$ for Hamiltonian identity, by utilizing the flexibility due to $f$.

\subsection{Notion of weak solutions and main result}
\label{sec-result}

%{\color{red}continue here...}

\begin{Definition} 
A function $\theta\in L^2_{loc}([0,\infty); \dot H^{-\frac12}(\mathbb T^2))$ is said to be a weak solution of (\ref{sqg}) with $f\in \dot H^{-r}$ if 
\begin{equation}\notag
\begin{split}
&\int_0^\infty\int_{\mathbb T^2} \Lambda^{-\frac12}\theta\partial_t\Lambda^{\frac12}\psi \, dxdt-\frac12\int_0^\infty\int_{\mathbb T^2}\Lambda^{-\frac12}\theta \Lambda^{\frac12}\left([\nabla^{\perp}\Lambda^{-1}, \nabla\psi]\theta \right)\, dxdt\\
=&\int_0^\infty\int_{\mathbb T^2}\Lambda^{-\frac12}\theta\Lambda^{\gamma+\frac12}\psi\, dxdt-\int_0^\infty\int_{\mathbb T^2}\Lambda^{-r} f\Lambda^r\psi\, dxdt
\end{split}
\end{equation}
for any $\psi\in C_0^\infty([0,\infty)\times \mathbb T^2)$ such that $\nabla\cdot \psi=0$. 
\end{Definition}
%&\int_0^\infty \left<\mathcal R_i^{\perp} \theta, \partial_t\Lambda^{-1}\phi^i\right>+\left<\mathcal R_j^{\perp}\theta, \mathcal R_i^{\perp}\Lambda^{-1}\theta\partial_j\phi^i\right>-\frac12\left< \mathcal R_i\mathcal R_j^{\perp}\theta, [\Lambda, \phi^i]\mathcal R_j^{\perp}\Lambda^{-1}\theta\right>\, dt\\
%=&\ \nu\int_0^\infty \left<\mathcal R_i^{\perp} \theta, \Lambda^{\gamma-1}\phi^i\right>\, dt-\int_0^\infty \left<\mathcal R_i^{\perp} \Lambda^{-\frac32} f, \Lambda^{\frac12}\phi^i\right>\, dt
Such weak solution in $L^2_{loc}([0,\infty); \dot H^{-\frac12}(\mathbb T^2))$ is well-defined thanks to 
\begin{equation}\notag
\int_{\mathbb T^2} \theta u\cdot \nabla\psi \, dx
=-\frac12\int_{\mathbb T^2} \theta [\nabla^{\perp} \Lambda^{-1}, \nabla\psi]\theta\, dx
\end{equation}
which is valid for $\theta\in\dot H^{-\frac12}$, see \cite{BSV}.

%Existence of weak solutions of (\ref{sqg}) in $\dot H^{-\frac12}$ without external forcing was established by Marchand \cite{Mar}. 

The main result is as follows. 
\begin{Theorem}\label{thm}
Let $0\leq \gamma<1-\alpha$, $-\frac12\leq \alpha<0$ and $0<\zeta<\frac12$. There exists $f\in C_t^0C_x^{2\alpha-1}$ such that there are more than one weak solutions $\theta \in C^0_tC_x^{\alpha}\cap C_t^{\zeta}C_x^0$ of (\ref{sqg}) with the external forcing $f$.
\end{Theorem}

%We mention that the weak solutions constructed here for the forced SQG have higher regularity ($C^{0-\varepsilon}$) than the solutions constructed in \cite{BSV, IM} for the SQG not driven by any forcing. The latter ones have regularity $C^{-\frac15-\varepsilon}$.

\medskip
In order to prove Theorem \ref{thm}, we first describe a general convex integration scheme for forced equation adapted to (\ref{sqg}) at the level $\Lambda^{-1}\theta$ in Section \ref{sec-outline}. Heuristic analysis is provided in Section \ref{sec-heuristics} for the SQG and forced SQG to give a brief outlook on how the regularity of the constructed solutions get improved with the presence of external forcing. In Section \ref{sec-prep} we lay out technical preparations and the main iterative proposition. 
Section \ref{sec-prop} is devoted to the proof of the main iterative argument. The proof of Theorem \ref{thm} is concluded in Section \ref{sec-proof}. 
%Finally we give an outline of an alternative proof of non-uniqueness for the stationary SQG in Appendix. 

\bigskip

\section{A convex integration scheme for forced equation}
\label{sec-outline}

%Denote $\eta=\Lambda^{-1}\theta$, hence $\theta=\Lambda\eta$, $u=\nabla^{\perp}\eta$. 
Constructions in this paper will be performed at the level of $\Lambda^{-1}\theta$, slightly different from the work of \cite{BSV} and \cite{IM}, where constructions occur at the level of $\Lambda^{-1}u$ (vector field) and $\theta$ respectively. In fact, a convex integration scheme at the level of $\Lambda^{-1}\theta$ was already applied in \cite{CKL} to construct non-trivial stationary solutions for SQG, where the error forcing term in the iterative approximating process appears in gradient form, see also \cite{DP} for the forced stationary SQG. However, as a byproduct of the present work, we show that the error forcing has a natural divergence form in the construction at the level of 
$\Lambda^{-1}\theta$. Therefore the scheme presented in this paper can provide a different proof for the result of \cite{CKL}. We leave the details for interested reader. 

We provide a sketch of a general iteration scheme for the forced evolutionary SQG based on our previous work \cite{DP} for the stationary forced SQG. 
It is not surprising that the presence of an external forcing relaxes the SQG equation. The authors of \cite{FK} exploited such flexibility to construct non-unique Leray-Hopf solutions for a forced dyadic model for the NSE, see also the work \cite{DF} for non-uniqueness construction for forced dyadic models for magnetohydrodynamics. Allowing an external forcing term,
one can find an initial pair $(\theta, u, f_1)$ and $(\widetilde\theta, \widetilde u, f_2)$ both satisfying (\ref{sqg}) with 
\[u=T[\theta], \ \ \widetilde u=T[\widetilde\theta], \ \ \theta\not\equiv \widetilde \theta\]
but not necessarily having $f_1\equiv f_2$. 
%Denote the sum $p=\frac12(\theta+\widetilde\theta)$ and the difference $m=\frac12(\theta-\widetilde\theta)$, and hence $\theta=p+m$ and $\widetilde\theta=p-m$. 
The sum $p=\frac12(\theta+\widetilde\theta)$ and the difference $m=\frac12(\theta-\widetilde\theta)$ satisfy the forced system
\begin{equation}\label{pm2}
\begin{split}
p_t+T[p]\cdot\nabla p+T[m]\cdot\nabla m=&-\nu \Lambda^\gamma p+\frac12(f_1+f_2),\\
m_t+T[p]\cdot\nabla m+T[m]\cdot\nabla p=&-\nu\Lambda^\gamma m+\frac12(f_1-f_2),\\
\nabla\cdot T[p]= 0, \ \ \nabla\cdot T[m]=&\ 0.
\end{split}
\end{equation} 
If $f_1\equiv f_2=f$, $(\theta, u, f)$ and $(\widetilde\theta, \widetilde u, f)$ are two distinct solutions of (\ref{sqg}) with the same forcing function. If $f_1\not\equiv f_2$, the forcing $\frac12(f_1-f_2)$ as an initial error in (\ref{pm2}) can be iteratively reduced by applying a convex integration scheme. %Three remarks are unfolded regarding the convex integration scheme in this context: (i) the convex integration scheme will be only applied to the $m$ equation not the $p$ equation in (\ref{pm2}); (ii) inspired by the work \cite{IM}, we recast the forcing $\frac12(f_1+f_2)=\nabla\cdot\nabla\cdot R$ and $\frac12(f_1-f_2)=\nabla\cdot\nabla\cdot \widetilde R$ into second derivative form; (iii) to improve the error estimates, we perform a special two-step construction in each iteration stage. Further details will be provided in Section \ref{sec-outline}.

Denote
\[\eta=\Lambda^{-1}\theta, \ \ \widetilde\eta=\Lambda^{-1}\widetilde\theta,\ \ \Pi=\frac12(\eta+\widetilde \eta), \ \ \mu=\frac12(\eta-\widetilde \eta),\]
and hence
\[p=\Lambda \Pi, \ \ m=\Lambda \mu, \ \ T[p]=\nabla^{\perp}\Pi, \ \ T[m]=\nabla^{\perp}\mu.\]
%Following the idea of \cite{CKL}, the convex integration is performed at the level of $\mu$. 
We rewrite (\ref{pm2}) as
\begin{equation}\label{pm3}
\begin{split}
\partial_t \Lambda \Pi+\nabla\cdot(\Lambda \Pi\nabla^{\perp}\Pi)+\nabla\cdot(\Lambda \mu\nabla^{\perp}\mu)=&- \Lambda^{\gamma+1} \Pi+\nabla\cdot\nabla\cdot \widetilde R,\\
\partial_t \Lambda \mu+\nabla\cdot(\Lambda \mu\nabla^{\perp}\Pi)+\nabla\cdot(\Lambda \Pi\nabla^{\perp}\mu)=&-\Lambda^{\gamma+1} \mu+\nabla\cdot\nabla\cdot   R\\
%\nabla\cdot \nabla^{\perp}\Pi= 0, \ \ \nabla\cdot \nabla^{\perp}\mu=&\ 0.
\end{split}
\end{equation} 
with stress tensors $\widetilde R$ and $ R$ satisfying 
\[\nabla\cdot\nabla\cdot \widetilde R=\frac12(f_1+f_2), \ \ \nabla\cdot\nabla\cdot  R=\frac12(f_1-f_2).\]
A convex integration scheme will be applied to construct a sequence of approximating solutions $\{(\Pi_q,\mu_q, \widetilde R_q, R_q)\}_{q\geq 0}$ of (\ref{pm3}) such that  $R_q\to 0$ which converge to $(\Pi, \mu, \widetilde R, 0)$ in a suitable sense as $q\to\infty$, with $\mu\not\equiv 0$. Thus it indicates the existence of two distinct  solutions $\theta=\Lambda(\Pi+\mu)$ and $\widetilde\theta=\Lambda(\Pi-\mu)$ of (\ref{sqg}) with forcing $f=\nabla\cdot\nabla\cdot \widetilde R$.

Let $q\geq 0$ be even and $(\Pi_q,\mu_q,  \widetilde R_q,  R_q)$ be the solution of (\ref{pm3}) at the $q$-th iteration, i.e.
\begin{equation}\label{pm-q}
\begin{split}
\partial_t \Lambda \Pi_q+\nabla\cdot(\Lambda \Pi_q\nabla^{\perp}\Pi_q)+\nabla\cdot(\Lambda \mu_q\nabla^{\perp}\mu_q)=&- \Lambda^{\gamma+1} \Pi_q+\nabla\cdot\nabla\cdot \widetilde R_q,\\
\partial_t \Lambda \mu_q+\nabla\cdot(\Lambda \mu_q\nabla^{\perp}\Pi_q)+\nabla\cdot(\Lambda \Pi_q\nabla^{\perp}\mu_q)=&-\Lambda^{\gamma+1} \mu_q+\nabla\cdot\nabla\cdot R_q
%\nabla\cdot \nabla^{\perp}\Pi_q= 0, \ \ \nabla\cdot \nabla^{\perp}\mu_q=&\ 0.
\end{split}
\end{equation} 
with the conversion
\[\Pi_q=\frac12\Lambda^{-1}(\theta_q+\widetilde\theta_q), \ \ \mu_q=\frac12\Lambda^{-1}(\theta_q-\widetilde\theta_q),\]
\[\theta_q=\Lambda (\Pi_q+\mu_q), \ \ \widetilde\theta_q=\Lambda (\Pi_q-\mu_q), \] 
\[ f_{1,q}=\nabla\cdot\nabla\cdot (\widetilde R_q+ R_q), \ \ f_{2,q}=\nabla\cdot\nabla\cdot (\widetilde R_q-R_q). \]
Note $(\theta_q, f_{1,q})$ and $(\widetilde\theta_q, f_{2,q})$ satisfy (\ref{sqg}). 
 
Taking advantage of the flexibility of the coupled system, each stage of the convex integration scheme consists two steps, from $q$-th to $(q+1)$-th step and from $(q+1)$-th to $(q+2)$-th step.  
In particular, we construct $W_{q+1}$ and $W_{q+2}$ to produce $(\Pi_{q+1}, \mu_{q+1})$ and $(\Pi_{q+2}, \mu_{q+2})$ respectively as 
\begin{equation}\label{two-step}
\begin{split}
\mu_{q+1}=&\ \mu_q+W_{q+1}, \ \ \Pi_{q+1}=\Pi_q-W_{q+1},\\
\mu_{q+2}=&\ \mu_{q+1}+W_{q+2}, \ \ \Pi_{q+2}=\Pi_{q+1}+W_{q+2}.
\end{split}
\end{equation}
Hence we have
\begin{equation}\notag
\begin{split}
\Lambda^{-1}\theta_{q+1}=&\ \eta_{q+1}= \Pi_{q+1}+\mu_{q+1}=\Pi_q+\mu_q\\
=&\ \eta_q=\Lambda^{-1}\theta_{q}, \\
\Lambda^{-1}\widetilde\theta_{q+1}=&\ \widetilde\eta_{q+1}= \Pi_{q+1}-\mu_{q+1}=\Pi_q-\mu_q-2W_{q+1}\\
=&\ \widetilde\eta_q-2W_{q+1}=\Lambda^{-1}\widetilde\theta_{q}-2W_{q+1},\\
\Lambda^{-1}\theta_{q+2}=&\ \eta_{q+2}= \Pi_{q+2}+\mu_{q+2}=\Pi_{q+1}+\mu_{q+1}+2W_{q+2}\\
=&\ \eta_{q+1}+2W_{q+2}=\Lambda^{-1}\theta_{q+1}+2W_{q+2}, \\
\Lambda^{-1}\widetilde\theta_{q+2}=&\ \widetilde\eta_{q+2}= \Pi_{q+2}-\mu_{q+2}=\Pi_{q+1}-\mu_{q+1}\\
=&\ \widetilde\eta_{q+1}=\Lambda^{-1}\widetilde\theta_{q+1}.
\end{split}
\end{equation}
The tuplet $(\Pi_{q+1}, \mu_{q+1}, \widetilde R_{q+1},  R_{q+1})$ satisfies (\ref{pm-q}) with $q$ replaced by $q+1$ for new stress tensors $\widetilde R_{q+1}$ and $R_{q+1}$, so does $(\Pi_{q+2}, \mu_{q+2}, \widetilde R_{q+2}, R_{q+2})$ for some $\widetilde R_{q+2}$ and $R_{q+2}$.
Subtraction of (\ref{pm-q}) at the $q$-th level and $(q+1)$-th level leads to
\begin{equation}\label{g-q1}
\begin{split}
\nabla\cdot\nabla\cdot R_{q+1}=&\ \partial_t\Lambda W_{q+1}+\nabla\cdot(\nabla^{\perp} \widetilde\eta_q \Lambda W_{q+1})+\nabla\cdot(\nabla^{\perp} W_{q+1}\Lambda \widetilde\eta_q)\\
&+\Lambda^{\gamma+1}W_{q+1}+\nabla\cdot\left(\nabla\cdot R_q-2\nabla^{\perp}W_{q+1}\Lambda W_{q+1}\right),\\
\nabla\cdot\nabla\cdot \widetilde R_{q+1}=&-\partial_t\Lambda W_{q+1}-\nabla\cdot(\nabla^{\perp} \widetilde\eta_q \Lambda W_{q+1})-\nabla\cdot(\nabla^{\perp} W_{q+1}\Lambda \widetilde\eta_q)\\
&-\Lambda^{\gamma+1}W_{q+1}
+\nabla\cdot(\nabla\cdot \widetilde R_q+2\nabla^{\perp}W_{q+1}\Lambda W_{q+1}),\\
\end{split}
\end{equation}
\begin{equation}\label{g-q2}
\begin{split}
\nabla\cdot\nabla\cdot R_{q+2}=&\ \partial_t\Lambda W_{q+2}+\nabla\cdot(\nabla^{\perp} \eta_{q+1} \Lambda W_{q+2})+\nabla\cdot(\nabla^{\perp} W_{q+2}\Lambda \eta_{q+1})\\
&+\Lambda^{\gamma+1}W_{q+2}+\nabla\cdot\left(\nabla\cdot R_{q+1}+2\nabla^{\perp}W_{q+2}\Lambda W_{q+2}\right),\\
\nabla\cdot\nabla\cdot \widetilde R_{q+2}=&\ \partial_t\Lambda W_{q+2}+\nabla\cdot(\nabla^{\perp} \eta_{q+1} \Lambda W_{q+2})+\nabla\cdot(\nabla^{\perp} W_{q+2}\Lambda \eta_{q+1})\\
&+\Lambda^{\gamma+1} W_{q+2}
+\nabla\cdot(\nabla\cdot \widetilde R_{q+1}+2\nabla^{\perp}W_{q+2}\Lambda W_{q+2}),
\end{split}
\end{equation}
According to (\ref{g-q1}) and (\ref{g-q2}), $W_{q+1}$ and $W_{q+2}$ will be constructed such that 
\begin{equation}\label{cancel1}
\begin{split}
\nabla\cdot R_q-2\nabla^{\perp}W_{q+1}\Lambda W_{q+1}\sim&\ \nabla^{\perp} F_{q+1},\\
\nabla\cdot R_{q+1}+2\nabla^{\perp}W_{q+2}\Lambda W_{q+2}\sim&\ \nabla^{\perp} F_{q+2}
\end{split}
\end{equation}
with small errors compared to other terms on the right hand side of (\ref{g-q1}) and some functions $F_{q+1}$ and $F_{q+2}$. 

Iterating the process for even $q\geq 2$ gives a sequence 
\[\{(\Pi_{q+i}, \mu_{q+i}, \widetilde R_{q+i}, R_{q+i})\}_{i=1,2; q\geq 0}\] 
with $(\Pi_{q+i}, \mu_{q+i}, \widetilde R_{q+i}, R_{q+i})$ satisfying the forced system (\ref{pm-q}). In particular, the force functions $R_{q+1}$ and
$R_{q+2}$ satisfy (\ref{g-q1}) and (\ref{g-q2}) respectively. 
The crucial point in this iteration process is
\begin{equation}\notag
\eta_{q+1}=\eta_q, \ \ \widetilde \eta_{q+2}=\widetilde \eta_{q+1}, \ \ \mbox{for any even} \ \ q\geq 0
\end{equation}
and consequently, since $u_{q}=T[\Lambda \eta_{q}]$, $\theta_q=\Lambda\eta_q$, $\widetilde u_{q}=T[\Lambda \widetilde \eta_{q}]$, $\widetilde \theta_q=\Lambda\widetilde\eta_q$, 
\begin{equation}\label{skip}
u_{q+1}=u_q, \ \ \widetilde u_{q+2}=\widetilde u_{q+1}, \ \ \theta_{q+1}=\theta_q, \ \ \widetilde \theta_{q+2}=\widetilde \theta_{q+1}, \ \ \mbox{for any even} \ \ q\geq 0
\end{equation}
which has the advantage to improve the transport errors and Nash errors in (\ref{g-q1}) and (\ref{g-q2}). 
%and hence improve the regularity of the constructed weak solutions. 

From (\ref{g-q1}) and (\ref{g-q2}) we observe that $ R_q$ is reduced to $ R_{q+1}$ by an amount which is gained in the process 
$\widetilde R_{q}\to \widetilde R_{q+1}$; 
%the ``reduced'' amount of force from $\widetilde R_q$ to $\widetilde R_{q+1}$ is gained by $R_{q+1}$ from $R_{q}$, 
while both $ R_{q+1}$ and $\widetilde R_{q+1}$ are ``reduced'' to  $R_{q+2}$ and $\widetilde R_{q+2}$ respectively by the same amount.  

%We mention that this convex integration scheme gives the same improvement for Nash error estimate as the alternating scheme used in \cite{BHP} for the forced Euler equation.

\bigskip

\section{Heuristics for SQG and forced SQG}
\label{sec-heuristics}

\subsection{Heuristic analysis for SQG}

For SQG, the authors of \cite{BSV} and \cite{IM} showed that
\begin{Theorem}\label{thm-sqg}
Let $0<\alpha<\frac{3}{10}$ and $0<\beta<\frac{1}{4}$. There exists nontrivial weak solutions $\theta$ to (unforced) SQG with compact support in time such that
\[\Lambda^{-\frac12}\theta\in C^0_tC_x^\alpha\cap C^\beta_t C^0_x.\]
\end{Theorem}
We provide a heuristic argument for the result. 
Consider the approximating system
\begin{equation}\label{eq-q}
\begin{split}
\partial_t\theta_q+\nabla\cdot(u_q\theta_q)=&\ \nabla\cdot\nabla\cdot R_q,\\
u_q=&\ T[\theta_q].
\end{split}
\end{equation}
%The stress vector $\widetilde R_q$ can be decomposed as 
%\[\widetilde R_q=c_{1,q}A_1+c_{2,q}A_2=: c_{1,q}A_1+R_q.\]
%Without loss of generality, we assume $|c_{1,q}|\geq |c_{2,q}|$. The goal is to construct a new solution such that the principal part $c_{1,q}A_1$ in the stress error gets reduced.

%We specify the index $I=(k, \pm)\in \mathbb Z\times \{\pm\}:= \Omega$. Denote $\bar I=(k, \mp)$. 
Let $\Omega$ be an index set to be specified in Section \ref{sec-prep}.
For $I\in\Omega$, 
%we choose phase functions $\xi_I$ satisfying $\xi_{\bar I}=-\xi_I$ and amplitude functions $a_{I,q+1}$ satisfying $a_{I,q+1}=a_{\bar I,q+1}$.
let $\theta_{I,q+1}$ and $\xi_{I}$ be the amplitude and phase functions respectively, satisfying
\[\theta_{\bar I}=\bar \theta_I, \ \ \xi_{\bar I}=-\xi_I.\]  
The phase functions $\xi_I$ are advected by $u_q$ with initial state $\xi_{I, in}$ on a short time interval with scale $\tau_q$, which is to be determined to optimize the error estimates. 
The increment $\Theta_{q+1}=\theta_{q+1}-\theta_q$ is constructed to take the form
\begin{equation}\notag
\begin{split}
\Theta_{I,q+1}=&\ e^{i\lambda\xi_{I}}(\theta_{I, q+1}+\delta\theta_{I,q+1}),\\
\Theta_{q+1}=&\sum_{I\in\Omega} \Theta_{I,q+1}
\end{split}
\end{equation}
where $\delta\theta_{I,q+1}$ is an error term.

Applying the Microlocal lemma, we have
\begin{equation}\notag
\begin{split}
U_{I,q+1}=&\ T[\Theta_{I,q+1}]= e^{i\lambda\xi_I}(u_{I,q+1}+\delta u_{I}), \ \ \
u_{I,q+1}= m(\nabla\xi_{I})\theta_{I,q+1}\\
U_{q+1}=&\ T[\Theta_{q+1}]=\sum_{I\in\Omega}T[\Theta_{I,q+1}]
\end{split}
\end{equation}
with $m(\xi)=-i\xi^{\perp}/|\xi|$ being the Fourier symbol for the operator $T=\nabla^{\perp}\Lambda^{-1}$.
Note 
\[u_{q+1}=u_q+U_{q+1}=T[\theta_q]+ T[\Theta_{q+1}].\]
The pair $(\theta_{q+1}, u_{q+1})$ is a solution of the next level approximating equation
%\begin{equation}\notag
%\partial_t\theta_{q+1}+\nabla\cdot(u_{q+1}\theta_{q+1})= \nabla\cdot\nabla\cdot R_{q+1}
%\end{equation}
with a new stress error $R_{q+1}$ satisfying
\begin{equation}\label{R-q1}
\begin{split}
\nabla\cdot\nabla\cdot R_{q+1}=& \left(\partial_t+u_{q}\cdot\nabla \right)\Theta_{q+1}+\nabla\cdot(U_{q+1}\theta_q)\\
&+ \nabla\cdot\sum_{J\neq \bar I}U_{J,q+1}\Theta_{I, q+1}\\
&+ \nabla\cdot\nabla\cdot\sum_{I\in\Omega}\left(U_{I,q+1}\Lambda^{-1}\Theta_{\bar I, q+1}+R_q\right)\\
=&: \nabla\cdot\nabla\cdot R_T+\nabla\cdot\nabla\cdot R_N+\nabla\cdot\nabla\cdot R_H+\nabla\cdot\nabla\cdot R_S.
\end{split}
\end{equation}
For parameters $\lambda_0\gg 1$, $b>1$ and $0<\beta<1$, define
\[\lambda_q=\left\lceil \lambda_0^{b^q}\right\rceil, \ \ \ q\in \mathbb N\cup \{0\}\]
and let $\delta_q=\lambda_q^{-\beta}$.
Assume $\theta_q$ and $u_q$ are localized to frequency $\sim \lambda_q$. The cancellation $U_{I,q+1}\Lambda^{-1}\Theta_{\bar I, q+1}+R_q$ suggests the scaling $|R_{q}|\sim \lambda_{q+1}^{-1}|\theta_{q+1}|^2$.
We make the inductive assumptions for some $L\geq 2$:
\begin{equation}\label{est-q-1}
\|\nabla^k u_q\|_{C^0}+\|\nabla^k \theta_q\|_{C^0}\lesssim \lambda_q^{k+\frac12}\delta_{q-1}^{\frac12}, \ \ k=1, 2, ..., L,
\end{equation}
\begin{equation}\label{est-q-2}
\|\nabla^k (\partial_t+u_q\cdot\nabla)u_q\|_{C^0}\lesssim \lambda_q^{k+2}\delta_{q-1}, \ \ k=0, 1, 2, ..., L-1,
\end{equation}
\begin{equation}\label{est-q-5}
\|\nabla^k R_{q}\|_{C^0}\lesssim \lambda_q^k\delta_{q}, \ \ k=1, 2, ..., L,
\end{equation}
\begin{equation}\label{est-q-6}
\|\nabla^k (\partial_t+u_q\cdot\nabla)R_{q}\|_{C^0}\lesssim \lambda_q^{k+\frac32}\delta_{q-1}^{\frac12}\delta_{q}, \ \ k=0, 1, 2, ..., L-1.
\end{equation}
Since 
\begin{equation}\notag
\begin{split}
R_T=&\ \div^{-2} P_{\sim \lambda_{q+1}} [(\partial_t+u_q\cdot\nabla)\Theta_{q+1}]\\
R_N=&\ \div^{-2} (U_{q+1}\cdot\nabla \theta_q),
\end{split}
\end{equation}
we expect
\begin{equation}\notag
\begin{split}
\|R_T\|_{C^0}\lesssim&\ \lambda_{q+1}^{-2} \|(\partial_t+u_q\cdot\nabla)\Theta_{q+1}\|_{C^0}\lesssim  \lambda_{q+1}^{-2}\tau_{q}^{-1}\lambda_{q+1}^{\frac12}\delta_q^{\frac12}=\lambda_{q+1}^{-\frac32}\tau_{q}^{-1}\delta_q^{\frac12},\\
\|R_N\|_{C^0}\lesssim&\ \lambda_{q+1}^{-2} \| U_{q+1}\cdot\nabla \theta_q\|_{C^0}\lesssim  \lambda_{q+1}^{-2}\lambda_{q+1}^{\frac12}\delta_q^{\frac12}\lambda_q\lambda_q^{\frac12}\delta_{q-1}^{\frac12}=\lambda_{q+1}^{-\frac32}\lambda_q^{\frac32}\delta_q^{\frac12}\delta_{q-1}^{\frac12}.
\end{split}
\end{equation}
On the other hand, one also expects to have
\begin{equation}\notag
\begin{split}
\|R_H\|_{C^0}=&\left\| \div^{-1}\sum_{J\neq \bar I}U_{J,q+1}\Theta_{I, q+1}\right\|_{C^0}\\
\lesssim&\ \lambda_{q+1}^{-1}\sum_{I}\| \theta_{I, q+1}\|_{C^0}^2\left(m(\|\nabla\xi_I)-m(\nabla\xi_{I,in})\|_{C^0}+\|\nabla\xi_I-\nabla \xi_{I,in}\|_{C^0} \right)\\
\lesssim&\ \lambda_{q+1}^{-1}\sum_{I}\| \theta_{I, q+1}\|_{C^0}^2\|\nabla\xi_I-\nabla \xi_{I,in}\|_{C^0} \\
\lesssim&\ \lambda_{q+1}^{-1}\sum_{I}\| \theta_{I, q+1}\|_{C^0}^2\lambda_q\tau_q \|u_q\|_{C^0} \\
\lesssim&\ \lambda_{q+1}^{-1}(\lambda_{q+1}^{\frac12}\delta_q^{\frac12})^2\lambda_q\tau_{q}\lambda_q^{\frac12}\delta_{q-1}^{\frac12}\\
=&\ \lambda_q^{\frac32}\tau_q\delta_q\delta_{q-1}^{\frac12}.
\end{split}
\end{equation}
To balance the error $R_T$ and $R_H$, we choose $\tau_{q}=\delta_{q-1}^{-\frac14}\delta_q^{-\frac14}\lambda_q^{-\frac34}\lambda_{q+1}^{-\frac34}$ such that
\[R_T\sim R_H\sim \lambda_{q+1}^{-\frac34}\lambda_q^{\frac34}\delta_q^{\frac34}\delta_{q-1}^{\frac14}.\]
In the end, we note $R_O$ is an error from the Microlocal lemma, which is supposed to be small compared to other error terms.  
Combining the estimates above gives
\begin{equation}\notag
\begin{split}
\|R_{q+1}\|_{C^0}\lesssim&\ \lambda_{q+1}^{-\frac34}\lambda_q^{\frac34}\delta_q^{\frac34}\delta_{q-1}^{\frac14}+\lambda_{q+1}^{-\frac32}\lambda_q^{\frac32}\delta_q^{\frac12}\delta_{q-1}^{\frac12}\\
\lesssim &\ \lambda_q^{-\frac34b+\frac34-\frac34\beta-\frac{\beta}{4b}}+\lambda_q^{-\frac32b+\frac32-\frac12\beta-\frac{\beta}{2b}}.
\end{split}
\end{equation}
To carry on the iterative process we need to make sure $\|R_{q+1}\|_{C^0}\lesssim \delta_{q+1}$, and hence require
\begin{equation}\label{para-conditions}
\begin{cases}
-\frac34b+\frac34-\frac34\beta-\frac{\beta}{4b}<-b\beta,\\
-\frac32b+\frac32-\frac12\beta-\frac{\beta}{2b}<-b\beta.
\end{cases}
\end{equation}
Thus we solve the first inequality in (\ref{para-conditions}), by recalling $b>1$
\begin{equation}\notag
\begin{split}
&-\frac34b+\frac34-\frac34\beta-\frac{\beta}{4b}+b\beta<0\\
\Longleftrightarrow&\ -\frac34(b-1)+\frac34\beta(b-1)+\frac{\beta}{4b}(b^2-1)<0\\
\Longleftrightarrow&\ -\frac34+\frac34\beta+\frac{\beta}{4b}(b+1)<0\\
\Longleftrightarrow&\ \beta<\frac{3b}{4b+1}.
\end{split}
\end{equation}
When $b=1^+$, we have $\beta<\frac35$. 
Similarly, the second inequality in (\ref{para-conditions}) is valid if
\begin{equation}\notag
\begin{split}
&-\frac32b+\frac32-\frac12\beta-\frac{\beta}{2b}+b\beta<0\\
\Longleftrightarrow&\ -\frac32(b-1)+\frac12\beta(b-1)+\frac{\beta}{2b}(b^2-1)<0\\
\Longleftrightarrow&\ -\frac32+\frac12\beta+\frac{\beta}{2b}(b+1)<0\\
\Longleftrightarrow&\ \beta<\frac{3}{1+\frac{b+1}{b}}
\end{split}
\end{equation}
which indicates $\beta<1$ for $b=1^+$.  The $C^\alpha$ regularity of $\Theta_{q+1}$ requires
\begin{equation}\notag
\|\Theta_{q+1}\|_{C^\alpha}\lesssim \lambda_{q+1}^\alpha \|\Theta_{q+1}\|_{C^0}\lesssim \lambda_{q+1}^{\alpha+\frac12}\delta_q^{\frac12}\lesssim \lambda_q^{b(\alpha+\frac12)-\frac12\beta}\lesssim 1
\end{equation}
which implies $\alpha<\frac{\beta}{2b}-\frac12<-\frac15$.

\subsection{Heuristic analysis for forced SQG}
\label{sec-heuristic-f}

In the case of forced SQG, thanks to (\ref{skip}), we have 
\begin{equation}\notag
u_{q+1}=u_q, \ \ \widetilde u_{q+2}=\widetilde u_{q+1}, \ \ \theta_{q+1}=\theta_q, \ \ \widetilde \theta_{q+2}=\widetilde \theta_{q+1}, \ \ \mbox{for any even} \ \ q\geq 0.
\end{equation}
Hence stress errors would be essentially
\begin{equation}\notag
\begin{split}
R_T=&\ \div^{-2} P_{\sim \lambda_{q+1}} [(\partial_t+u_{q-1}\cdot\nabla)\Theta_{q+1}]\\
R_N=&\ \div^{-2} (U_{q+1}\cdot\nabla \theta_{q-1})\\
\end{split}
\end{equation}
and 
\begin{equation}\notag
\|R_H\|_{C^0}
\lesssim \lambda_{q+1}^{-1}\sum_{I}\| \theta_{I, q+1}\|_{C^0}^2\lambda_{q-1}\tau_{q-1} \|u_{q-1}\|_{C^0}.
\end{equation}
Consequently we have
\begin{equation}\notag
\begin{split}
\|R_T\|_{C^0}\lesssim&\ \lambda_{q+1}^{-2} \|(\partial_t+u_{q-1}\cdot\nabla)\Theta_{q+1}\|_{C^0}\lesssim  \lambda_{q+1}^{-2}\tau_{q-1}^{-1}\lambda_{q+1}^{\frac12}\delta_q^{\frac12}=\lambda_{q+1}^{-\frac32}\tau_{q-1}^{-1}\delta_q^{\frac12},\\
\|R_N\|_{C^0}\lesssim&\ \lambda_{q+1}^{-2} \| U_{q+1}\cdot\nabla \theta_{q-1}\|_{C^0}\lesssim  \lambda_{q+1}^{-2}\lambda_{q+1}^{\frac12}\delta_q^{\frac12}\lambda_{q-1}\lambda_{q-1}^{\frac12}\delta_{q-2}^{\frac12}=\lambda_{q+1}^{-\frac32}\lambda_{q-1}^{\frac32}\delta_q^{\frac12}\delta_{q-2}^{\frac12},\\
\|R_H\|_{C^0}
%\lesssim&\ \lambda_{q+1}^{-1}\sum_{I}\| \theta_{I, q+1}\|_{C^0}^2\lambda_q\tau_q \|u_q\|_{C^0} \\
\lesssim&\ \lambda_{q+1}^{-1}(\lambda_{q+1}^{\frac12}\delta_q^{\frac12})^2\lambda_{q-1}\tau_{q-1}\lambda_{q-1}^{\frac12}\delta_{q-2}^{\frac12}\lesssim
 \lambda_{q-1}^{\frac32}\tau_{q-1}\delta_q\delta_{q-2}^{\frac12}.
\end{split}
\end{equation}
By requiring $\|R_T\|_{C^0}\sim \|R_H\|_{C^0}$ we have 
\[\tau_{q-1}= \delta_{q-2}^{-\frac14}\delta_q^{-\frac14}\lambda_{q-1}^{-\frac34}\lambda_{q+1}^{-\frac34}  \]
and hence 
\[\|R_T\|_{C^0}\sim \|R_H\|_{C^0}\lesssim \lambda_{q+1}^{-\frac32} \delta_q^{\frac12}  \delta_{q-2}^{\frac14}\delta_q^{\frac14}\lambda_{q-1}^{\frac34}\lambda_{q+1}^{\frac34}= \lambda_{q+1}^{-\frac34}\lambda_{q-1}^{\frac34}\delta_q^{\frac34}\delta_{q-2}^{\frac14}.\]
It follows that
\begin{equation}\notag
\begin{split}
\|R_{q+1}\|_{C^0}\leq&\ \|R_T\|_{C^0}+\|R_N\|_{C^0}+\|R_H\|_{C^0}+\|R_O\|_{C^0}\\
\lesssim&\ \lambda_{q+1}^{-\frac34}\lambda_{q-1}^{\frac34}\delta_q^{\frac34}\delta_{q-2}^{\frac14}+\lambda_{q+1}^{-\frac32}\lambda_{q-1}^{\frac32}\delta_q^{\frac12}\delta_{q-2}^{\frac12}\\
\lesssim&\ \lambda_q^{-\frac34b+\frac{3}{4b}-\frac34\beta-\frac{1}{4b^2}\beta}+\lambda_q^{-\frac32b+\frac3{2b}-\frac12\beta-\frac{1}{2b^2}\beta}.
\end{split}
\end{equation}
Imposing $\|R_{q+1}\|_{C^0}\leq \delta_{q+1}=\lambda_q^{-b\beta}$ leads to 
\begin{equation}\notag
\begin{split}
&\begin{cases}
-\frac34b+\frac{3}{4b}-\frac34\beta-\frac{1}{4b^2}\beta+b\beta<0,\\
-\frac32b+\frac3{2b}-\frac12\beta-\frac{1}{2b^2}\beta+b\beta<0
\end{cases}\\
\Longleftrightarrow &
\begin{cases}
-\frac{3}{4b}(b^2-1)+\frac{3}{4}\beta(b-1)+\frac{1}{4b^2}\beta(b^3-1)<0,\\
-\frac{3}{2b}(b^2-1)+\frac{1}{2}\beta(b-1)+\frac{1}{2b^2}\beta(b^3-1)<0
\end{cases}\\
\Longleftrightarrow &
\begin{cases}
-\frac{3}{4b}(b+1)+\frac{3}{4}\beta+\frac{1}{4b^2}\beta(b^2+b+1)<0,\\
-\frac{3}{2b}(b+1)+\frac{1}{2}\beta+\frac{1}{2b^2}\beta(b^2+b+1)<0
\end{cases}\\
\Longleftrightarrow &
\begin{cases}
\beta<\frac{3b(b+1)}{4b^2+b+1},\\
\beta<\frac{3b(b+1)}{2b^2+b+1}.
\end{cases}\\
\end{split}
\end{equation}
Choosing $b$ close enough to $1^+$ gives $\beta<1$. Hence $\theta$ is expected to be in $C^\alpha$ with $\alpha<\frac{\beta}{2b}-\frac12<0$.

\bigskip

\section{Technical preparation}
\label{sec-prep}

\subsection{The index set}
Adapting notation from \cite{IM}, let $F=\{\pm(1,2), \pm(2,1)\}$ and $\mathbb F=F/(+,-)=\{(1,2), (2,1)\}$.
The index $I$ takes the form $I=(k,v)\in \mathbb Z\times F$.  Define a wave $W_I=W_{(k,v)}$ with phase function $\xi_I$ and amplitude function $a_I$, where $k$ indicates the interval of time support and $v$ the oscillation direction. We also denote the conjugate index $\bar I=(k,\bar v)\in \mathbb Z\times F$ and choose conjugate pairs of waves $\{W_I, W_{\bar I}\}$ such that $W_{\bar I}=\bar W_I$, $\xi_{\bar I}=-\xi_I$ and $a_{\bar I}=\bar a_I$. Hence the summation in (\ref{def-w}) is real-valued. Moreover, denote
\begin{equation}\notag
[k]=
\begin{cases}
0 \ \ \mbox{if} \ \ k \ \ \mbox{mod} \ 2=0\\
1 \ \ \mbox{if} \ \ k \ \ \mbox{mod} \ 2=1.
\end{cases}
\end{equation}

%\medskip

\subsection{Flow map}

For any $j\in\mathbb Z$, we define $\Phi_j(x,t)$ as the solution to the flow map
\begin{equation}\label{transport-phi}
\begin{split}
(\partial_t+u_q\cdot\nabla) \Phi_j=&\ 0,\\
\Phi_j(x, j\tau_{q})=&\ x.
\end{split}
\end{equation}
We also define $G_{q,j}$ as the solution to the transport equation
\begin{equation}\label{transport-g}
\begin{split}
(\partial_t+u_q\cdot\nabla) G_{q,j}=&\ 0,\\
G_{q,j}(x,j\tau_{q})=&\ G_q(x,j\tau_{q}).
\end{split}
\end{equation}
The following transport estimates can be found in \cite{BSV}.
\begin{Lemma}\label{le-phi}
Let $\Phi_j$ be a solution to (\ref{transport-phi}). It satisfies
\begin{equation}\notag
\begin{split}
\|D\Phi_j(t)-Id\|_{C^0}\leq&\ e^{(t-t_0)\|Du_q\|_{C^0}}-1\leq (t-t_0)\|Du_q\|_{C^0}e^{(t-t_0)\|Du_q\|_{C^0}},\\
\|D^N\Phi_j(t)\|_{C^0}\leq&\ C (t-t_0)\|D^Nu_q\|_{C^0}e^{C(t-t_0)\|Du_q\|_{C^0}}, \ \ N\geq 2.
\end{split}
\end{equation}
\end{Lemma}

%\medskip

\subsection{Microlocal Lemma}

We recall a microlocal lemma from \cite{IV}.
\begin{Lemma}\label{le-micro}
Let $\lambda\in \mathbb Z$, $\xi:\mathbb T^2\to \mathbb C$ and $a:\mathbb T^2\to \mathbb C$ be smooth functions. Let $K:\mathbb R^2\to \mathbb C$ be a Schwartz function. Define the operator
\begin{equation}\notag
T[W](x)=\int_{\mathbb R^2}W(x-h)K(h)\, dh, \ \ W(x)=e^{i\lambda\xi(x)}a(x).
\end{equation}
Then we have
\begin{equation}\notag
T[W](x)=e^{i\lambda\xi(x)}\left(a(x)\widehat K(\lambda\nabla \xi(x))+\delta[TW](x) \right)
\end{equation}
with the error term given by
\begin{equation}\notag
\begin{split}
\delta[TW](x)=&\int_0^1 dr\frac{d}{dr} \int_{\mathbb R^2}e^{-i\lambda\nabla\xi(x)\cdot h}e^{iZ(r,x,h)}a(x-rh)K(h) \, dh\\
Z(r,x,h)=&\ r\lambda\int_0^1 h^lh^{l'}\nabla_l\nabla_{l'}\xi(x-sh)(1-s)\, ds.
\end{split}
\end{equation}
\end{Lemma}

In analogy with the bilinear microlocal lemma of \cite{IM}, we establish the bilinear microlocal lemma at a different level for the quadratic drift term.

\begin{Lemma} \label{le-bilinear}
Let $W_{I}$ and $W_{\bar I}$ be conjugate waves defined as
\[W_I(x,t)=P_{\lambda}\left(a_I(x,t)e^{i\lambda \xi_{I}(x,t)}\right), \ \ W_{\bar I}(x,t)=P_{\lambda}\left(a_{I}(x,t)e^{i\lambda \xi_{\bar I}(x,t)}\right), \ \ \xi_I=-\xi_{\bar I}\]
with $P_{\lambda}$ being the frequency localization operator defined as in Subsection \ref{sec-frequency}.
Denote 
\begin{equation}\notag
Q(x,t)= \nabla^{\perp}W_I\Lambda W_{\bar I}+\nabla^{\perp}W_{\bar I}\Lambda W_{I}.
\end{equation}
Then we have 
\[Q^l=\nabla_j B^{jl}_{\lambda}[W_I, W_{\bar I}] \]
for a bilinear operator $B_\lambda$ in
the expansion form
\begin{equation}\notag
B_{\lambda}[W_I, W_{\bar I}](x,t)=a_I^2(x,t)\widehat K_{\lambda,sym}(\lambda\nabla\xi_I, -\lambda\nabla\xi_I)+\delta B_I(x,t)
\end{equation}
with $\widehat K_{\lambda,sym}$ and the error term $\delta B_I(x,t)$ to be described below.
\end{Lemma}
\pf
For general
\begin{equation}\notag
Q(x,t)= \nabla^{\perp}W_I\Lambda W_{J}+\nabla^{\perp}W_{J}\Lambda W_{I},
\end{equation}
the Fourier transform of $Q$ is given by
\begin{equation}\notag
\begin{split}
\widehat Q(\xi)=&\ i\int_{\mathbb R^2}(\xi-k)^{\perp}\widehat W_I(\xi-k)|k|\widehat W_J(k)\, dk+i\int_{\mathbb R^2}|\xi-k|\widehat W_I(\xi-k)k^{\perp}\widehat W_J(k)\, dk\\
=&\ i\int_{\mathbb R^2}\left((\xi-k)^{\perp}|k|+k^{\perp}|\xi-k|\right)_{\approx \lambda}\widehat W_I(\xi-k)\widehat W_J(k)\, dk\\
=&\ i\int_{\mathbb R^2}|\xi-k||k|\left(\frac{(\xi-k)^{\perp}}{|\xi-k|}+\frac{k^{\perp}}{|k|}\right)_{\approx \lambda}\widehat W_I(\xi-k)\widehat W_J(k)\, dk\\
=&\int_{\mathbb R^2}|\xi-k||k|\left(m(\xi-k)+m(k)\right)_{\approx \lambda}\widehat W_I(\xi-k)\widehat W_J(k)\, dk\\
=&\int_{\mathbb R^2}\left((\xi-k)^{\perp}|k|-(-k)^{\perp}|\xi-k|\right)_{\approx \lambda}\widehat W_I(\xi-k)\widehat W_J(k)\, dk\\
\end{split}
\end{equation}
with $m(\xi)=i\xi^{\perp}/|\xi|$. Since $m$ is odd, it follows
\begin{equation}\notag
\begin{split}
\widehat Q(\xi)
=&\int_{\mathbb R^2}|\xi-k||k|\left(m(\xi-k)-m(-k)\right)_{\approx \lambda}\widehat W_I(\xi-k)\widehat W_J(k)\, dk\\
=&\int_{\mathbb R^2}\xi_j  \int_0^1 \nabla^j m_{\approx \lambda}^l(u_\sigma)\, d\sigma |\xi-k||k|\widehat W_I(\xi-k)\widehat W_J(k)\, dk\\
\end{split}
\end{equation}
with $u_\sigma=\sigma(\xi-k)-(1-\sigma) k$. Therefore we have
\begin{equation}\notag
\widehat Q^l=i\xi_j\int_{\mathbb R^2}\widehat K_{\lambda}^{jl}(\zeta,k)\widehat W_I(\zeta)\widehat W_J(k)\, dk
\end{equation}
with $\zeta=\xi-k$ and 
\begin{equation}\label{kernel-def}
\widehat K_{\lambda}^{jl}(\zeta, k)=-i|\zeta||k|\int_0^1\nabla^j m_{\approx \lambda}^l(u_\sigma)\widehat \chi_{\approx \lambda}(\zeta)\widehat\chi_{\approx \lambda}(k)\, d\sigma.
\end{equation}
Denote the symmetric part of $\widehat K_{\lambda}^{jl}$ as 
\[\widehat K_{\lambda, sym}^{jl}=\frac12(\widehat K_{\lambda}^{jl}+\widehat K_{\lambda}^{lj})\]
and 
\begin{equation}\notag
\widehat Q^l_s=i\xi_j\int_{\mathbb R^2}\widehat K_{\lambda,sym}^{jl}(\zeta,k)\widehat W_I(\zeta)\widehat W_J(k)\, dk.
\end{equation}
Denote $K_{\lambda,sym}(h_1,h_2)$ by the physical space kernel as the inverse Fourier transform of $\widehat K_{\lambda,sym}(\zeta,k)$, that is
\begin{equation}\notag
\widehat K_{\lambda,sym}(\zeta,k)=\int_{\mathbb T^2\times \mathbb T^2} e^{-i(\zeta, k)\cdot (h_1,h_2)}K_{\lambda,sym}(h_1,h_2)\, dh_1dh_2.
\end{equation}
Taking the inverse Fourier transform on $\widehat Q_s$ gives
\begin{equation}\notag
Q_s^l=\nabla_j\int_{\mathbb T^2\times \mathbb T^2}K^{jl}_{\lambda, sym}(h_1,h_2) W_I(x-h_1)W_J(x-h_2)\, dh_1dh_2.
\end{equation}
Thus we define the bilinear operator $B_{\lambda}[W_I, W_{J}]$ as 
\begin{equation}\notag
B^{jl}_{\lambda}[W_I, W_{J}]=\int_{\mathbb T^2\times \mathbb T^2}K^{jl}_{\lambda, sym}(h_1,h_2) W_I(x-h_1)W_J(x-h_2)\, dh_1dh_2.
\end{equation}
From (\ref{kernel-def}) we have, noting $\widehat \chi_{\approx \lambda} (\pm\lambda\nabla\xi_I)=1$
\begin{equation}\notag
\begin{split}
&\widehat K_{\lambda}^{jl}(\lambda\nabla \xi_I, -\lambda\nabla \xi_I)\\
=&-i \lambda^2|\nabla\xi_I|^2\int_0^1\nabla^jm^l_{\approx \lambda}(\lambda\nabla\xi_I)\widehat \chi_{\approx \lambda} (\lambda\nabla\xi_I) \widehat \chi_{\approx \lambda} (-\lambda\nabla\xi_I) \, d\sigma\\
=&-i \lambda^2|\nabla\xi_I|^2 \nabla^jm^l_{\approx \lambda}(\lambda\nabla\xi_I).
\end{split}
\end{equation}
Direct computation yields 
\begin{equation}\label{kernel-s}
\begin{split}
&\widehat K_{\lambda, sym}^{jl}(\lambda\nabla \xi_I, -\lambda\nabla \xi_I)\\
=&-i \frac12\lambda^2|\nabla\xi_I|^2 \left( \nabla^jm^l +\nabla^lm^j\right)(\lambda\nabla\xi_I)\\
=&\ \frac{1}{2}  \lambda |\nabla \xi_I|^{-1} 
\begin{pmatrix}
-2\partial_1\xi_I\partial_2\xi_I & (\partial_1\xi_I)^2 - (\partial_2\xi_I)^2 \\
(\partial_1\xi_I)^2 - (\partial_2\xi_I)^2 & 2\partial_1\xi_I\partial_2\xi_I
\end{pmatrix}.
\end{split}
\end{equation}
% where we used $\xi_I^{\perp}=i \xi_I$.  
We observe that $\widehat K_{\lambda, sym}$ is homogeneous of degree 1 in both $\lambda$ and $\nabla\xi_I$.  
Consequently for 
\[W_I(x,t)=P_{\lambda}\left(a_I(x,t)e^{i\lambda \xi_{I}(x,t)}\right), \ \ W_{\bar I}(x,t)=P_{\lambda}\left(a_{I}(x,t)e^{i\lambda \xi_{\bar I}(x,t)}\right)\]
we deduce, by using $\xi_{\bar I}=-\xi_{I}$
\begin{equation}\notag
\begin{split}
B_{\lambda}[W_I, W_{\bar I}]=&\int_{\mathbb T^2\times \mathbb T^2}e^{i\lambda\xi_I(x-h_1)} e^{-i\lambda\xi_{I}(x-h_2)}a_I(x-h_1)a_{I}(x-h_2)K_{\lambda, sym}(h_1,h_2)\, dh_1dh_2\\
=&\int_{\mathbb T^2\times \mathbb T^2}e^{i\lambda(\xi_I(x-h_1)-\xi_I(x))} e^{-i\lambda (\xi_I(x-h_2)-\xi_I(x))}\\
&\cdot a_I(x-h_1)a_{I}(x-h_2)K_{\lambda, sym}(h_1,h_2)\, dh_1dh_2.
\end{split}
\end{equation}
In view of Taylor expansion we have
\begin{equation}\notag
\begin{split}
\xi_I(x-h)-\xi_I(x)=&-\nabla\xi_I(x)h+\int_{0}^1h_jh_l\partial_j\partial_l\xi_I(x-sh)(1-s)\, ds,\\
a_I(x-h)-a_I(x)=&-\nabla a_I(x)\cdot h+\int_{0}^1h_jh_l\partial_j\partial_l a_I(x-sh)(1-s)\, ds.
\end{split}
\end{equation}
Note that $|h_1|, |h_2|\sim \lambda^{-1}$ and $\lambda\gg 1$. Hence
\[a_I(x-h_1)\sim a_I(x), \ \ a_I(x-h_2)\sim a_I(x).\]
Therefore the leading order term of $B_{\lambda}[W_I, W_{\bar I}]$ is 
\begin{equation}\notag
\begin{split}
&\int_{\mathbb T^2\times \mathbb T^2}e^{i\lambda\nabla\xi_I h_1} e^{-i\lambda\nabla\xi_{I}h_2}a_I^2(x)K_{\lambda, sym}(h_1,h_2)\, dh_1dh_2\\
%=& \ a_I^2(x)\widehat K(\lambda\xi_I \nabla\Phi(x), \lambda\xi_{\bar I} \nabla\Phi(x))\\
=& \ a_I^2(x)\widehat K_{\lambda, sym}(\lambda \nabla\xi_I , -\lambda\nabla\xi_{I}).
\end{split}
\end{equation}
More precisely we can write 
\begin{equation}\notag
\begin{split}
B_{\lambda}[W_I, W_{\bar I}]=& \ a_I^2(x)\widehat K_{\lambda, sym}(\lambda\nabla\xi_I, -\lambda\nabla\xi_{I})\\
&+\int_{\mathbb T^2\times \mathbb T^2}e^{-i\lambda\nabla\xi_Ih_1}e^{i\lambda\nabla\xi_Ih_2} K_{\lambda, sym}(h_1,h_2) a_I(x) Y(x,h_1)\, dh_1dh_2\\
&+\int_{\mathbb T^2\times \mathbb T^2}e^{-i\lambda\nabla\xi_Ih_1}e^{i\lambda\nabla\xi_Ih_2} K_{\lambda, sym}(h_1,h_2) a_I(x) Y(x,h_2)\, dh_1dh_2\\
&+\int_{\mathbb T^2\times \mathbb T^2}e^{-i\lambda\nabla\xi_Ih_1}e^{i\lambda\nabla\xi_Ih_2} K_{\lambda, sym}(h_1,h_2) Y(x,h_1)Y(x,h_2)\, dh_1dh_2\\
=&: B_{\lambda, p}[W_I, W_{\bar I}]+B_{\lambda,1}[W_I, W_{\bar I}]+B_{\lambda,2}[W_I, W_{\bar I}]+B_{\lambda,3}[W_I, W_{\bar I}]
\end{split}
\end{equation}
with 
\begin{equation}\notag
Y(x,h)=\int_0^1\frac{d}{dr}e^{ir\lambda\int_0^1h^jh^l\partial_j\partial_l\xi_I(x-sh)(1-s)\, ds}a_I(x-rh)\, dr.
\end{equation}
Thus we have
\begin{equation}\notag
B_{\lambda}[W_I, W_{\bar I}](x,t)=a_I^2(x,t)\widehat K_{\lambda,sym}(\lambda\nabla\xi_I, -\lambda\nabla\xi_I)+\delta B_I(x,t)
\end{equation}
with $\widehat K_{\lambda,sym}(\lambda\nabla\xi_I, -\lambda\nabla\xi_I)$ given by (\ref{kernel-s}) and 
\[\delta B_I(x,t)= B_{\lambda,1}[W_I, W_{\bar I}]+B_{\lambda,2}[W_I, W_{\bar I}]+B_{\lambda,3}[W_I, W_{\bar I}].\]

\cbdu

%{\color{blue}This bilinear microlocal lemma can be used to give a different proof of the non-uniqueness of the stationary SQG, in the form 
%\[\nabla^\perp \eta \Lambda \eta+\nabla^\perp p=\nabla \cdot R.\]}

%\medskip

\subsection{Second order anti-divergence operator}
\label{sec-anti-div}

For smooth and mean-zero function $g:\mathbb C^2\to \mathbb C$, the second order anti-divergence operator is defined as
\begin{equation}\notag
\mathcal R=2\mathcal R_1-\mathcal R_2
\end{equation}
with 
\[\mathcal R_1^{jl}=\nabla^j\nabla^l\Delta^{-2}, \ \ \ \mathcal R_2^{jl}=\delta^{jl}\Delta^{-1}.\]
It is easy to observe that $\mathcal R[g]$ is a symmetric traceless tensor and $\nabla_j\nabla_l \mathcal R^{jl}[g]=g$.

\medskip

\subsection{Main iteration step}
\label{sec-iteration}

%Fix a large constant $\lambda_0>0$. Let $b>1$. Choose the frequency number as the integer
%\[\lambda_q=\left\lceil \lambda_0^{b^q}\right\rceil, \ \ \ q\in \mathbb N\cup \{0\}.\]
%The magnitude measure is given by $\delta_q=\lambda_q^{-\beta}$ for a parameter $\beta>0$ to be specified later. 
Let $\lambda_q$ and $\delta_q$ be defined as in Section \ref{sec-heuristics}, denoting the frequency number and magnitude measure respectively.
%We also choose the frequency localization number $r_{q+1}=(\lambda_q\lambda_{q+1})^{\frac12}$. 
In the process of iteration described in Subsection \ref{sec-outline}, we set the scale  
\[|R_q|\sim \delta_q.\]
The increment $W_{q+1}$ will be constructed such that (i) it is supported in Fourier space near the frequency $\lambda_{q+1}$; (ii) it has $C^{\alpha}$ regularity for some $\alpha$ to be determined. 
Thus in view of the cancelation in (\ref{cancel1})
\[\nabla\cdot  R_q-2\nabla^{\perp} W_{q+1}\Lambda W_{q+1}\]
which would essentially lead to the scaling
\[ |R_q| \sim \lambda_{q+1}^{-1}|\nabla^{\perp} W_{q+1}\Lambda W_{q+1}|\]
by invoking the bilinear microlocal Lemma \ref{le-bilinear}.
Thus we expect to have
\begin{equation}\label{w-size}
|W_{q+1}|\sim \left(\lambda_{q+1}^{-1}\delta_{q}\right)^{\frac12}.
\end{equation}
The $C^{\alpha}$ regularity requirement for $W_{q+1}$ indicates 
\[\lambda_{q+1}^{\alpha} \left(\lambda_{q+1}^{-1}\delta_{q}\right)^{\frac12}\lesssim 1\]
which implies 
\begin{equation}\label{alpha}
\alpha<\frac12+\frac{\beta}{2b}.
\end{equation}
Once the appropriate value for $\beta$ is achieved for the inductive process to be carried on, (\ref{alpha}) determines the regularity of the constructed solutions.

To ease notation, we denote $\|\cdot\|_{C^N}=\|\cdot\|_{N}$ in the rest of the paper. 
For an integer $L\geq 1$ and constant $M>0$, we assume the inductive estimates for any $q\geq 0$
%\begin{equation}\label{induct-eta-0}
%\|\eta_q\|_{0}\leq M(1-\lambda_q^{-\frac12}\delta_{q-1}^{\frac12})
%\end{equation}
\begin{equation}\label{induct-theta}
\|\mu_q\|_{N+1}+\|\Pi_q\|_{N+1}
\leq M\lambda_q^{N+\frac12}\delta_{q-1}^{\frac12}, \ \ \ 0\leq N\leq L
\end{equation}
%\begin{equation}\label{induct-eta}
%\|\eta_q\|_{N}\leq M\lambda_q^{N-\frac12}\delta_{q-1}^{\frac12},  \ \ \ N\in\{1, 2, ..., L\}
%\end{equation}
\begin{equation}\label{induct-R}
\|R_q\|_{N}\leq \lambda_q^N\delta_q, \ \ \ 0\leq N\leq L
\end{equation}
\begin{equation}\label{induct-mu-t}
\begin{split}
&\|(\partial_t+\nabla^{\perp}\Pi_q\cdot\nabla)\Lambda\mu_q\|_{N}+\|(\partial_t+\nabla^{\perp}\Pi_q\cdot\nabla)\Lambda\Pi_q\|_{N}\\
&+\|(\partial_t+\nabla^{\perp}\mu_q\cdot\nabla)\Lambda\mu_q\|_{N}+\|(\partial_t+\nabla^{\perp}\mu_q\cdot\nabla)\Lambda\Pi_q\|_{N}\\
\leq&\ M\lambda_q^{N+2}\delta_{q-1}, \ \ \ 0\leq N\leq L-1
\end{split}
\end{equation}
\begin{equation}\label{induct-R-t}
\begin{split}
&\|(\partial_t+\nabla^{\perp}\Pi_q\cdot\nabla)R_q\|_{N}+\|(\partial_t+\nabla^{\perp}\mu_q\cdot\nabla)R_q\|_{N}\\
\leq&\ \lambda_q^{N+\frac32}\delta_q\delta_{q-1}^{\frac12}, \ \ \ 0\leq N\leq L-1
\end{split}
\end{equation}
and 
\begin{equation}\label{induct-freq}
\supp \widehat \mu_q\cup \supp \widehat \Pi_q\subset \{\xi: |\xi|\leq \lambda_q\}.
\end{equation}
In view of the convention
\[\theta_q=\Lambda\eta_q=\Lambda(\Pi_q+\mu_q), \ \ u_q=\nabla^{\perp}\eta_q=\nabla^{\perp}(\Pi_q+\mu_q),\]
\[\widetilde\theta_q=\Lambda\widetilde\eta_q=\Lambda(\Pi_q-\mu_q), \ \ \widetilde u_q=\nabla^{\perp}\widetilde \eta_q=\nabla^{\perp}(\Pi_q-\mu_q)\]
and the crucial feature (\ref{skip}) of the scheme outlined in Section \ref{sec-outline}, it follows from the inductive estimates above that for even $q$
\begin{equation}\label{induct-u1}
\|\widetilde u_q\|_{N}+\|\widetilde \theta_q\|_{N}
\leq M\lambda_{q-1}^{N+\frac12}\delta_{q-2}^{\frac12}, \ \ \ 0\leq N\leq L
\end{equation}
\begin{equation}\label{induct-u2}
\|u_{q+1}\|_{N}+\|\theta_{q+1}\|_{N}
\leq M\lambda_q^{N+\frac12}\delta_{q-1}^{\frac12}, \ \ \ 0\leq N\leq L
\end{equation}
\begin{equation}\label{induct-mat1}
\|(\partial_t+\widetilde u_q\cdot\nabla) R_q\|_{N}
\lesssim \lambda_q^N\lambda_{q-1}^{\frac32}\delta_{q-2}^{\frac12}\delta_q
\end{equation}
\begin{equation}\label{induct-mat2}
\|(\partial_t+u_{q+1}\cdot\nabla) R_{q+1}\|_{N}
\lesssim \lambda_{q+1}^N\lambda_q^{\frac32}\delta_{q-1}^{\frac12}\delta_{q+1}.
\end{equation}

The main inductive argument is given in the following proposition. 
\begin{Proposition}\label{prop-main}
Fix $L\geq 2$ and $M\geq 4$. For even $q$, let $(\Pi_q, \mu_q, \widetilde R_q, R_q)$ be a solution of (\ref{pm-q}) satisfying (\ref{induct-theta})-(\ref{induct-freq}) with $\supp_tR_q\subset J$ for a closed non-empty interval $J\subset \mathbb R$. There exist tuples $(\Pi_{q+1}, \mu_{q+1}, \widetilde R_{q+1}, R_{q+1})$ and $(\Pi_{q+2}, \mu_{q+2}, \widetilde R_{q+2}, R_{q+2})$ in the form of 
\begin{equation}\notag
\begin{split}
\mu_{q+1}=&\ \mu_q+W_{q+1}, \ \ \ \Pi_{q+1}=\Pi_q-W_{q+1}, \\
\mu_{q+2}=&\ \mu_{q+1}+W_{q+2}, \ \ \Pi_{q+2}=\Pi_{q+1}+W_{q+2}
\end{split}
\end{equation}
%with \[\theta_{q+1}=\Lambda \eta_{q+1}=\Lambda \eta_{q}+\Lambda W_{q+1}, \ \ u_{q+1}=\nabla^{\perp}\eta_{q+1}=\nabla^{\perp}\eta_{q}+\nabla^{\perp}W_{q+1}\]
satisfying system (\ref{pm-q}) and estimates (\ref{induct-theta})-(\ref{induct-R-t}) with $q$ replaced by $q+1$ and $q+2$ respectively. Moreover, we have
\begin{equation}\label{time-supp}
\supp_t R_{q+1}\cup \supp_t W_{q+1} \subset \left\{ t+h: t\in J, |h|\leq 4\lambda_q^{-\frac32}\delta_{q-1}^{-\frac12}\right\}.
\end{equation}

\end{Proposition}

\bigskip

\section{Proof of Proposition \ref{prop-main}}
\label{sec-prop}

It is sufficient to show the iteration estimates in Proposition \ref{prop-main} from the $q$-th level to $(q+1)$-th level; the iteration from the $(q+1)$-th level to $(q+2)$-th level can be established analogously. 

\subsection{Building blocks for the increment}
\label{sec-blocks}

We consider the increment $W_{q+1}$ in the form
\begin{equation}\label{def-w}
W_{q+1}(x,t)=\sum_{I}P_{I,\lambda_{q+1}}\left(a_{I,q+1}(x,t)e^{i\lambda_{q+1}\xi_I(x,t)} \right)
\end{equation}
with the frequency localization operator $P_{I, \lambda_{q+1}}$ and amplitude function $a_{I,q+1}$ to be specified later. The phase function $\xi_I$ with $I=(k,v)$ is the solution of the backward flow equation 
\begin{equation}\notag
\begin{cases}
\partial_t \xi_I+u_q\cdot\nabla \xi_I=0,\\
\xi_I(x,k\tau_q)=\xi_{I, in}(x)
\end{cases}
\end{equation}
where the initial time is $k\tau_q$ with the time step $\tau_q$ to be determined later, and the initial data is 
\[\xi_{I, in}(x)=J^k f\cdot x \ \ \ \mbox{with}  \ J \ \mbox{being a rotation of} \ 90^\circ .\]
In particular, the time step $\tau_q$ will be chosen sufficiently small such that 
\begin{equation}\label{time-condition}
\left|\nabla \xi_I-\nabla\xi_{I, in} \right|\leq c_1, \ \ \ |\nabla\xi_I|\geq 1
\end{equation}
for a small constant $c_1$. 

Denote 
\begin{equation}\notag
\begin{split}
D_{t,q}=&\ \partial_t+u_q\cdot\nabla, \ \ \widetilde D_{t,q}=\partial_t+\widetilde u_q\cdot\nabla, \\
D_{t,q,\Pi}=&\ \partial_t+\nabla^{\perp}\Pi_q\cdot\nabla, \ \ D_{t,q,\mu}=\partial_t+\nabla^{\perp}\mu_q\cdot\nabla.
\end{split}
\end{equation}
\begin{Lemma}\label{le-phase} \cite{IM}
For an appropriate choice of the time scale $\tau_q$, condition (\ref{time-condition}) holds for $|t-k\tau_q|\leq \tau_q$. Moreover, the estimates
\begin{equation}\notag
\begin{split}
\|\nabla_a D^r_{t,q}(\nabla\xi_I)\|_{0}\lesssim& \left( \lambda_{q+1}\lambda_q^{-1}\right)^{\frac{(|a|+1-L)_+}{L}}\lambda_q^{|a|} \left( \lambda_q^{\frac32}\delta_{q-1}^{\frac12}\right)^r, \ \ \ |a|\geq 0, \ \ r=0,1\\
\|\nabla_a D^2_{t,q}(\nabla\xi_I)\|_{0}\lesssim& \left( \lambda_{q+1}\lambda_q^{-1}\right)^{\frac{(|a|+2-L)_+}{L}}\lambda_q^{|a|} \left( \lambda_q^{\frac32}\delta_{q-1}^{\frac12}\right)^2, \ \ \ |a|\geq 0
\end{split}
\end{equation}
hold. 
\end{Lemma}

\begin{Remark}
For $r=1$ the estimate is valid when the differential operator $\nabla_{a_1}D_{t,q}\nabla_{a_2}$ is applied to $\nabla\xi_I$ with $|a_1|+|a_2|=|a|$. Similarly, in the case of $r=2$, the estimate holds when applying $\nabla_{a_1}D_{t,q}\nabla_{a_2}D_{t,q}\nabla_{a_3}$ to $\nabla\xi_I$ with $|a_1|+|a_2|+|a_3|=|a|$.
\end{Remark}

For convenience we denote
\begin{equation}\label{def-wj}
\widetilde W_{I,q+1}=a_{I,q+1}(x,t)e^{i\lambda_{q+1}\xi_I(x,t)}, \ \   W_{I,q+1}=P_{I,\lambda_{q+1}}\widetilde W_{I, q+1}
\end{equation}
and hence
\begin{equation}\notag
W_{q+1}(x,t)=\sum_{I} W_{I, q+1}.
\end{equation}

\medskip

\subsection{Time cutoffs}
\label{sec-time}
Let $0\leq \chi\leq 1$ be a smooth cutoff function with $\chi\equiv 1$ on $[1,2]$ and $\chi\equiv 0$ on the complement of $[\frac12, 4]$, which gives a partition of unity in the sense of
\[\sum_{k\in \mathbb Z}\chi^2(t-k)=1, \ \ \ \forall \ \ t\in\mathbb R.\]
Denote the rescaled function
\begin{equation}\notag
\chi_{k,q}(t):=\chi_k(t)=\chi(t\tau_{q}^{-1}-k).
\end{equation}
For a given index $I=(k,v)$, the amplitude function $a_{I, q+1}$ will be constructed consisting the time cutoff function $\chi_{k,q}(t)$ such that $W_{I,q+1}$ will be supported on the time interval $[k\tau_q-\frac23\tau_q, k\tau_q+\frac23\tau_q]$.

\medskip

\subsection{Frequency localization}
\label{sec-frequency}
As conventionally in other convex integration work,  the operator $P_{I, \lambda_{q+1}}$ in (\ref{def-w}) is defined for a continuous function $g$ and index $I$
\begin{equation}\notag
P_{I, \lambda_{q+1}} g=\int_{\mathbb R^2} g(x-h)\phi_{I,\lambda_{q+1}}(h)\, dh
\end{equation}
with a smooth bump function $\phi_I\in C^{\infty}_c(B_{|\nabla\xi_{I, in}|/2}(\nabla\xi_{I,in}))$ satisfying 
\[\phi_I(\xi)=1 \ \ \ \mbox{if} \ \ |\nabla \xi_I-\nabla\xi_{I,in}|\leq \frac14 |\nabla\xi_{I,in}|\]
and $\phi_{I,\lambda_{q+1}}=\phi_I(\lambda_{q+1}^{-1}\xi)$. Therefore $W_{I,q+1}$ and $W_{q+1}$ are supported in frequency space near $\lambda_{q+1}$, that is
\[\supp \widehat W_{I,q+1}\cup \supp \widehat W_{q+1} \subset \{\xi\in \mathbb R^2: \frac12\lambda_{q+1}\leq |\xi| \leq 2\lambda_{q+1} \}. \]

\medskip

\subsection{Mollification of the stress tensor}
\label{sec-molli}

To avoid derivative loss, we mollify the stress tensor in space and in time along the coarse flow. 

Let $l_q$ be a space length scale to be determined later. 
Let $\rho: \mathbb R^2\to \mathbb R$ be a Schwartz kernel and its rescale 
\[\rho_{l_q}(h)=l_q^{-2}\rho(l_q^{-1}h)\]
satisfying 
\[ \int_{\mathbb R^2} \rho(h)\, dh=1\]
and the vanishing moment condition
\[\int_{\mathbb R^2} h^{a} \rho_{l_q}(h)\, dh=0, \ \ \mbox{for all multi-index} \ a \ \mbox{with} \ 1\leq |a| \leq L. \]
We first define the spatial mollification of stress tensor $R_q$ as 
\[R_{q, l_q}=\rho_{l_q}* \rho_{l_q}* R_q.\]

Let $\Phi_s(x,t)$ be the flow map of $\partial_t+u_q\cdot \nabla$, that is, $\Phi_s(x,t)=(\Phi_s(x,t), t+s)$ is the solution to 
\begin{equation}\notag
\begin{cases}
\frac{d}{ds}\Phi_s(x,t)= u_q(\Phi_s(x,t)),\\
\Phi_0(x,t)= (x,t).
\end{cases}
\end{equation}
Let $\tilde \rho$ be a standard temporal mollifier compactly supported on $[-1,1]$ with $\int_{\mathbb R}\tilde\rho(s)\,ds=1$. Take $\tilde \rho_{\tilde \tau_{q}}(s)=\tilde \tau_{q}^{-1}\tilde\rho(\tilde\tau_{q}^{-1}s)$ for a time length scale $\tilde\tau_q\leq (\lambda_q^{\frac32}\delta_{q-1}^{\frac12})^{-1}$.

We then define the mollified stress tensor $\bar R_q$ as
\begin{equation}\notag
\bar R_q=\int_{-\tilde\tau_q}^{\tilde\tau_q} R_{q,l_q}(\Phi_s(x,t))\tilde\rho_{\tilde\tau_q}(s)\,ds.
\end{equation}
We choose the mollification scales 
\begin{equation}\label{scale-molli-xt}
l_q=c_0\left(\frac{\lambda_{q+1}}{\lambda_q}\right)^{-\frac{3}{2L}}\lambda_q^{-1}, \ \ \tilde\tau_q=c_0\left(\frac{\lambda_{q+1}}{\lambda_q}\right)^{-\frac{3}{2}}\lambda_q^{-\frac32}\delta_q^{-\frac12}
\end{equation}
for a constant $c_0$ depending only on $\rho$ and $\tilde \rho$.

\begin{Lemma}\label{le-molli}\cite{IM}
The estimates
\begin{equation}\notag
\begin{split}
\|\nabla_a D^r_{t,q, \Pi, \mu}\bar R_q\|_{0}\lesssim& \left( \lambda_{q+1}\lambda_q^{-1}\right)^{\frac{3(|a|+r-L)_+}{2L}}\lambda_q^{|a|} \left( \lambda_q^{\frac32}\delta_{q-1}^{\frac12}\right)^r\delta_q, \\
& \ \ \ \ \ \ \ \ \ \ \ \ \ \ \ \ \ \ \ \mbox{for}\ \ |a|\geq 0, \ \ r=0,1\\
\|\nabla_a D^2_{t,q, \Pi, \mu}\bar R_q\|_{0}\lesssim& \left( \lambda_{q+1}\lambda_q^{-1}\right)^{\frac{3(|a|+1-L)_+}{2L}}\lambda_q^{|a|} \left( \lambda_q^{\frac32}\delta_{q-1}^{\frac12}\right)\tilde \tau_q^{-1} \delta_q, \\
&\ \ \ \ \ \ \ \ \ \ \ \ \ \ \ \ \ \ \ \mbox{for}\ \ |a|\geq 0\\
\|R_q-\bar R_q\|_{0}\lesssim& \left(\left( \lambda_q^{\frac32}\delta_{q-1}^{\frac12}\right)\tilde \tau_q+l_q^L\lambda_q^L \right) \delta_q
\end{split}
\end{equation}
hold for the mollified stress tensor.
\end{Lemma}

\medskip

\subsection{Lifting function}
\label{sec-lift}
We define a lifting function $e(t)$ with the size of $R_q$ which is compactly supported in time, such as
\begin{equation}\notag
e^{\frac12}(t)=(C \delta_q)^{\frac12} \tilde\rho_{(\lambda_q^{\frac32}\delta_{q-1}^{\frac12})^{-1}}* \chi_{N_0(J)}(t)
\end{equation}
where $N_0(J)=\{t+h: t\in J, |h|\leq 2(\lambda_q^{\frac32}\delta_{q-1}^{\frac12})^{-1}\}$ and $C$ is a large constant. We note $e(t)\geq C\delta_q$ for $t\in \{s+h| s\in J, |h|\leq \tilde\tau_q+ (\lambda_q^{\frac32}\delta_{q-1}^{\frac12})^{-1}\}$ for an appropriate $\tilde\tau_q$. Moreover, the time scale $\tilde\tau_q$ is chosen such that 
\begin{equation}\notag
\left \|\left( \frac{d}{dt}e^{\frac12}(t)\right)^r \right\|_{0}\lesssim (\lambda_q^{\frac32}\delta_{q-1}^{\frac12})^r \delta_q^{\frac12}, \ \ 0\leq r\leq 2,
\end{equation}
and 
\begin{equation}\notag
\supp_t e(t)\subset \{t+h| t\in J, |h|\leq 3(\lambda_q^{\frac32}\delta_{q-1}^{\frac12})^{-1}\}.
\end{equation}
The lifting function $e(t)$ is designed to treat the low frequency stress error, see Section \ref{sec-amplitude}. 

\medskip

\subsection{Amplitude functions}
\label{sec-amplitude}

The amplitude functions $a_{I,q+1}$ will be chosen such that cancellation occurs in (\ref{cancel1}). In fact, we define
\begin{equation}\label{coeff-a}
a_{I,q+1}(x,t)=\lambda_{q+1}^{-\frac12}b_{I,q+1}(x,t)\chi_{k}(t)e^{\frac12}(t)
\end{equation}
for the lifting function $e(t)$, time cut-off functions $\chi_{k}(t)$, and the coefficient function $b_{I,q+1}(x,t)$ to be defined later. The purpose is to establish

\begin{Lemma}\label{le-decomp}
There exist a constant matrix $M_{[k]}$ and coefficient functions $b_{I,q+1}(x,t)$ such that
\begin{equation}\notag
\sum_{I=(k,v)\atop \in \mathbb Z\times F}a_{I,q+1}^2(x,t)\widehat {K}_{\lambda_{q+1}, sym}(\lambda_{q+1}\nabla\xi_I,-\lambda_{q+1}\nabla\xi_I)-e(t)\chi_k^2(t)\left(M_{[k]}-\frac{\bar R_q}{e(t)}\right)=0.
\end{equation}
\end{Lemma}
\pf
We first aim to identify the constant matrix $M_{[k]}$ by solving the linear system
\begin{equation}\label{matrix-linear}
\sum_{I=(k,v)\atop \in \{k\}\times F}a_{I,q+1}^2(x,t)\widehat {K}_{\lambda_{q+1}, sym}(\lambda_{q+1}\nabla\xi_{I,in},-\lambda_{q+1}\nabla\xi_{I,in})=e(t)\chi_k^2(t)M_{[k]}
\end{equation}
with constant coefficients $b_{I,q+1}(x,t)=\hat b_{I,q+1}$ (to be determined) and the initial conditions $\xi_{I,in}$ defined in Subsection \ref{sec-blocks}. It follows from (\ref{kernel-s}) and direct computation that
\begin{equation}\notag
\widehat K_{1,sym}(\nabla\xi_{(k, J^k(1,2), in)}, -\nabla\xi_{(k, J^k(1,2), in)})=
\begin{cases}
\frac{1}{2\sqrt 5}
\begin{bmatrix}
-4 & -3 \\
-3 & 4
\end{bmatrix} \ \ \mbox{if} \ \ [k]=0 \\
\frac{1}{2\sqrt 5}
\begin{bmatrix}
4 & 3 \\
3 & -4
\end{bmatrix} \ \ \mbox{if} \ \ [k]=1
\end{cases}
\end{equation}
\begin{equation}\notag
\widehat K_{1,sym}(\nabla\xi_{(k, J^k(2,1), in)}, -\nabla\xi_{(k, J^k(2,1), in)})=
\begin{cases}
\frac{1}{2\sqrt 5}
\begin{bmatrix}
-4 & 3 \\
3 & 4
\end{bmatrix} \ \ \mbox{if} \ \ [k]=0 \\
\frac{1}{2\sqrt 5}
\begin{bmatrix}
4 & -3 \\
-3 & -4
\end{bmatrix} \ \ \mbox{if} \ \ [k]=1.
\end{cases}
\end{equation}
Noting 
\[\widehat K_{\lambda_{q+1}, sym}(\lambda_{q+1}\nabla\xi_{I,in}, -\lambda_{q+1}\nabla\xi_{I,in})=\lambda_{q+1}^{-1}\widehat K_{1, sym}(\nabla\xi_{I,in}, -\nabla\xi_{I,in}), \]
we conclude that (\ref{matrix-linear}) is solved for constant coefficients $\hat b_{I, q+1}=1$ and the constant matrices $M_{[k]}$ given by
\begin{equation}\notag
M_{[k]}=\frac{1}{2\sqrt 5} 
\begin{bmatrix}
-8 & 0 \\
0 & 8
\end{bmatrix}
\ \ \mbox{for} \ \ k\equiv 0 \ \ \mbox{mod} \ \ 2,
\end{equation}
\begin{equation}\notag
M_{[k]}=\frac{1}{2\sqrt 5}
\begin{bmatrix}
8 & 0 \\
0 & -8
\end{bmatrix}
\ \ \mbox{for} \ \ k\equiv 1 \ \ \mbox{mod} \ \ 2.
\end{equation}

We proceed to solve the nonlinear equation in the lemma. 
Denote 
\begin{equation}\notag
\nabla_j\nabla_lR^{jl}_S=\nabla_j\nabla_l\bar R^{jl}_q+\nabla_l \sum_{I=(k,v)\in \mathbb Z\times F}\nabla^{\perp}_j W_{I,q+1} \Lambda W_{\bar I, q+1}
\end{equation}
It follows from the bilinear microlocal Lemma \ref{le-bilinear} that
\begin{equation}\notag
\begin{split}
\nabla_j\nabla_lR^{jl}_S=&\ \nabla_j\nabla_l \sum_{I=(k,v)\in \mathbb Z\times F}B^{jl}_{\lambda_{q+1}}[W_{I,q+1},W_{\bar I, q+1}] \\
&-\nabla_j\nabla_l\sum_{k\in \mathbb Z}\chi^2_k(t)\left(e(t)M_{k}^{jl}-\bar R_q^{jl}\right)
\end{split}
\end{equation}
and 
\begin{equation}\notag
\begin{split}
R_s^{jl}=& \sum_{I=(k,v)\atop \in \mathbb Z\times F}a_{I,q+1}^2(x,t)\widehat {K}^{jl}_{\lambda_{q+1}, sym}(\lambda_{q+1}\nabla\xi_I,-\lambda_{q+1}\nabla\xi_I)\\
&-\sum_{k\in\mathbb Z}e(t)\chi_k^2(t)\left(M_{[k]}-\frac{\bar R^{jl}_q}{e(t)}\right)+ \sum_{I=(k,v)\atop \in \mathbb Z\times F}\delta B_I^{jl}(x,t).
\end{split}
\end{equation}
We will show that there exist appropriate coefficient functions $b_{I, q+1}$ such that 
\begin{equation}\label{main-cancel}
\begin{split}
&\sum_{I=(k,v)\atop \in \mathbb Z\times F}a_{I,q+1}^2(x,t)\widehat {K}^{jl}_{\lambda_{q+1}, sym}(\lambda_{q+1}\nabla\xi_I,-\lambda_{q+1}\nabla\xi_I)\\
=&\sum_{k\in\mathbb Z}e(t)\chi_k^2(t)\left(M_{[k]}-\frac{\bar R^{jl}_q}{e(t)}\right)
\end{split}
\end{equation}
and hence
\[R_S^{jl}= \sum_{I=(k,v)\atop \in \mathbb Z\times F}\delta B_I^{jl}.\]
Note that equation (\ref{main-cancel}) holds if for any $k\in \mathbb Z$ we can solve
\begin{equation}\label{main-cancel-1}
\sum_{I\in\{k\}\times \mathbb F} b_{I,q+1}^2\widehat K_{1,sym}[\nabla\xi_I,-\nabla\xi_I]= M_{[k]}-\frac{\bar R_q}{e(t)}.
\end{equation}
It follows from similar analysis as in \cite{IM} that we can choose a small enough constant $c_1$ satisfying (\ref{time-condition}) and a large enough constant $K$ satisfying $e(t)>K\|\bar R_q\|_{C^0}$ so that equation (\ref{main-cancel-1}) is a small perturbation of the linear system (\ref{matrix-linear}), and hence (\ref{main-cancel-1}) can be solved.

\cbdu

\begin{Lemma}
For any $I\in \mathbb Z\times F$, the following estimates hold for $\frac{\bar R_q}{e(t)}$ and the coefficient functions $b_{I,q+1}$
\begin{equation}\notag
\|\nabla_{\vec a}D_{t,q,\Pi,\mu}^r\left(\frac{\bar R_q}{e(t)}\right)\|_{0}+\|\nabla_{\vec a}D_{t,q,\Pi,\mu}^rb_{I,q+1}\|_{0}\lesssim \left(\frac{\lambda_{q+1}}{\lambda_q}\right)^{\frac{3}{2L}(|\vec a|+1-L)_+}\lambda_q^{|\vec a|} \left( \lambda_q^{\frac32}\delta_{q-1}^{\frac12}\right)^r
\end{equation}
for $|\vec a|\geq 0$, $r=0,1$ and 
\begin{equation}\notag
\|\nabla_{\vec a}D_{t,q,\Pi,\mu}^2\left(\frac{\bar R_q}{e(t)}\right)\|_{0}+\|\nabla_{\vec a}D_{t,q,\Pi,\mu}^2b_{I,q+1}\|_{0}\lesssim \left(\frac{\lambda_{q+1}}{\lambda_q}\right)^{\frac{3}{2L}(|\vec a|+1-L)_+}\lambda_q^{|\vec a|} \left( \lambda_q^{\frac32}\delta_{q-1}^{\frac12}\right)\tilde \tau_q^{-1}
\end{equation}
for $|\vec a|\geq 0$.
\end{Lemma}
The proof follows from a similar analysis as in \cite{IM}.

\medskip

By the microlocal Lemma \ref{le-micro} we have 
\begin{equation}\notag
\begin{split}
W_{I,q+1}=&\ P_{I,\lambda_{q+1}}\left(a_{I,q+1}(x,t)e^{i\lambda_{q+1}\xi_I(x,t)} \right)=e^{i\lambda_{q+1}\xi_I(x,t)}(a_{I,q+1}+\delta a_{I,q+1})\\
\Lambda W_{I,q+1}=&\ e^{i\lambda_{q+1}\xi_I(x,t)}(\theta_{I,q+1}+\delta \theta_{I,q+1})\\
 \nabla^{\perp} W_{I,q+1} =&\ i\lambda_{q+1}\nabla^{\perp} \xi_I W_{I,q+1}+e^{i\lambda_{q+1}\xi_I(x,t)}(u_{I,q+1}+\delta u_{I,q+1})
\end{split}
\end{equation}
with
\[\theta_{I,q+1}=\lambda_{q+1}|\nabla \xi_I|a_{I,q+1}, \ \ \ u_{I,q+1}=\nabla^{\perp}a_{I,q+1}.\]

\begin{Lemma}\label{le-amplitude-main}
For any $I\in \mathbb Z\times F$, we have the following estimates for the principal amplitude functions
\begin{equation}\notag
\begin{split}
%\|\nabla_{\vec a}D_t^ra_{I,q+1}\|_{C^0}
%\lesssim&\ \lambda_{q+1}^{-\frac12}\left( \frac{\lambda_{q+1}}{\lambda_q}\right)^{\frac{3}{2L}(|\vec a|+1-L)_+}\lambda_q^{|\vec a|}\delta_q^{\frac12}\tau_q^{-r}, \ \ |\vec a|\geq 0, \ \ r=0,1 \\
%\|\nabla_{\vec a}D_t^2a_{I,q+1}\|_{C^0}
%\lesssim&\ \lambda_{q+1}^{-\frac12}\left( \frac{\lambda_{q+1}}{\lambda_q}\right)^{\frac{3}{2L}(|\vec a|+1-L)_+}\lambda_q^{|\vec a|}\delta_q^{\frac12}\left( \lambda_q^{\frac32}\delta_{q-1}^{\frac12}\right)\tilde \tau_q^{-1}, \ \ |\vec a|\geq 0 \\
&\lambda_{q+1}\|\nabla_{\vec a}D_{t,q, \Pi,\mu}^ra_{I,q+1}\|_{0}+\|\nabla_{\vec a}D_{t,q, \Pi,\mu}^r\theta_{I,q+1}\|_{0}+\|\nabla_{\vec a}D_{t,q, \Pi,\mu}^ru_{I,q+1}\|_{0}\\
\lesssim&\ \lambda_{q+1}^{\frac12}\left( \frac{\lambda_{q+1}}{\lambda_q}\right)^{\frac{3}{2L}(|\vec a|+1-L)_+}\lambda_q^{|\vec a|}\delta_q^{\frac12}\tau_q^{-r}, \ \ |\vec a|\geq 0, \ \ r=0,1 \\
&\lambda_{q+1}\|\nabla_{\vec a}D_{t,q, \Pi,\mu}^2a_{I,q+1}\|_{0}+\|\nabla_{\vec a}D_{t,q, \Pi,\mu}^2\theta_{I,q+1}\|_{0}+\|\nabla_{\vec a}D_{t,q, \Pi,\mu}^2u_{I,q+1}\|_{0}\\
\lesssim&\ \lambda_{q+1}^{\frac12}\left( \frac{\lambda_{q+1}}{\lambda_q}\right)^{\frac{3}{2L}(|\vec a|+1-L)_+}\lambda_q^{|\vec a|}\delta_q^{\frac12}\left( \lambda_q^{\frac32}\delta_{q-1}^{\frac12}\right)\tilde \tau_q^{-1}, \ \ |\vec a|\geq 0,
\end{split}
\end{equation}
\begin{equation}\notag
\begin{split}
&\lambda_{q+1}\|\nabla_{\vec a}\widetilde D_{t,q}^ra_{I,q+1}\|_{0}+\|\nabla_{\vec a}\widetilde D_{t,q}^r\theta_{I,q+1}\|_{0}+\|\nabla_{\vec a}\widetilde D_{t,q}^ru_{I,q+1}\|_{0}\\
\lesssim&\ \lambda_{q+1}^{\frac12}\left( \frac{\lambda_{q+1}}{\lambda_q}\right)^{\frac{3}{2L}(|\vec a|+1-L)_+}\lambda_q^{|\vec a|}\delta_q^{\frac12}\tau_{q-1}^{-r}, \ \ |\vec a|\geq 0, \ \ r=0,1 \\
&\lambda_{q+1}\|\nabla_{\vec a}\widetilde D_{t,q}^2a_{I,q+1}\|_{0}+\|\nabla_{\vec a}\widetilde D_{t,q}^2\theta_{I,q+1}\|_{0}+\|\nabla_{\vec a}\widetilde D_{t,q}^2u_{I,q+1}\|_{0}\\
\lesssim&\ \lambda_{q+1}^{\frac12}\left( \frac{\lambda_{q+1}}{\lambda_q}\right)^{\frac{3}{2L}(|\vec a|+1-L)_+}\lambda_q^{|\vec a|}\delta_q^{\frac12}\left( \lambda_{q-1}^{\frac32}\delta_{q-2}^{\frac12}\right)\tilde \tau_{q-1}^{-1}, \ \ |\vec a|\geq 0,
\end{split}
\end{equation}
\end{Lemma}
\pf
The estimates for $a_{I, q+1}$ can be obtained similarly as in the proof of Proposition 4.4 in \cite{IM}.
While the estimates of $\theta_{I, q+1}$ and $u_{I, q+1}$ follows from the fact
\[|\theta_{I, q+1}|\lesssim \lambda_{q+1}|a_{I, q+1}|, \ \ |u_{I, q+1}|\lesssim \lambda_{q+1}|a_{I, q+1}|.\]

\cbdu

\medskip

Estimates for the error terms generated by microlocal lemma:

\begin{Lemma}\label{le-error1-micro}
Let $L\geq 2$. Assume 
%the coarse-scale velocity $u_{q}$ satisfies 
\begin{equation}\notag%\label{coarse}
\lambda_q\|\nabla_{\vec a}\mu_q\|_{0}+\lambda_q\|\nabla_{\vec a}\Pi_q\|_{0}
\lesssim \left(\frac{\lambda_{q+1}}{\lambda_q} \right)^{\frac{1}{L}(|\vec a|-L)_+}\lambda_q^{|\vec a|}\lambda_q^{\frac12}\delta_{q-1}^{\frac12}, \ \ |\vec a|\geq 0. 
\end{equation}
For $\delta B_I$ defined in the bilinear microlocal lemma, we have for all $I\in\mathbb Z\times F$
\begin{equation}\notag
\begin{split}
\|\nabla_{\vec a}D_{t,q,\Pi,\mu}^{\rho} \delta B_I\|_{0}\lesssim&\ \lambda_{q+1}^{-2} \left(\frac{\lambda_{q+1}}{\lambda_q}\right)^{-1}\left(\frac{\lambda_{q+1}}{\lambda_q} \right)^{\frac{3}{2L}(|\vec a|+2-L)_+}\lambda_q^{|\vec a|}\delta_q\tau_q^{-\rho},\\
\|\nabla_{\vec a}\widetilde D_{t,q}^{\rho} \delta B_I\|_{0}\lesssim&\ \lambda_{q+1}^{-2} \left(\frac{\lambda_{q+1}}{\lambda_q}\right)^{-1}\left(\frac{\lambda_{q+1}}{\lambda_q} \right)^{\frac{3}{2L}(|\vec a|+2-L)_+}\lambda_q^{|\vec a|}\delta_q\tau_{q-1}^{-\rho}
\end{split} 
\end{equation}
for $ |\vec a|\geq 0$ and $0\leq \rho\leq 1$.
\end{Lemma}
We refer to \cite{IM} for a proof of the lemma, noting the $\delta B_I$ here is a multiple of $\lambda_{q+1}^{-2}$ of the error $\delta B_I$ from \cite{IM}.

\begin{Lemma}\label{le-error2-micro}
For $\delta a_{I,q+1}$, $\delta \theta_{I,q+1}$ and $\delta u_{I,q+1}$ appeared from the microlocal lemma, we have for all $I\in\mathbb Z\times F$
\begin{equation}\notag
\begin{split}
&\lambda_{q+1}\|\nabla_{\vec a}D_{t,q, \Pi,\mu}^r\delta a_{I,q+1}\|_{0}+\|\nabla_{\vec a}D_{t,q, \Pi,\mu}^r \delta \theta_{I,q+1}\|_{0}+\|\nabla_{\vec a}D_{t,q, \Pi,\mu}^r \delta u_{I,q+1}\|_{0}\\
\lesssim&\ \lambda_{q+1}^{\frac12} \left(\frac{\lambda_{q+1}}{\lambda_q}\right)^{-1}\left(\frac{\lambda_{q+1}}{\lambda_q} \right)^{\frac{3}{2L}(|\vec a|+2-L)_+}\lambda_q^{|\vec a|}\delta_q^{\frac12}\tau_q^{-r}, \ \ |\vec a|\geq 0, \ \ r=0,1 \\
&\lambda_{q+1}\|\nabla_{\vec a}D_{t,q, \Pi,\mu}^2\delta a_{I,q+1}\|_{0}+\|\nabla_{\vec a}D_{t,q, \Pi,\mu}^2 \delta \theta_{I,q+1}\|_{0}+\|\nabla_{\vec a}D_{t,q, \Pi,\mu}^2 \delta u_{I,q+1}\|_{0}\\
\lesssim&\ \lambda_{q+1}^{\frac12} \left(\frac{\lambda_{q+1}}{\lambda_q}\right)^{-1}\left(\frac{\lambda_{q+1}}{\lambda_q} \right)^{\frac{3}{2L}(|\vec a|+2-L)_+}\lambda_q^{|\vec a|}\delta_q^{\frac12}\left(\lambda_q^{\frac32}\delta_{q-1}^{\frac12}\right)\tilde \tau_q^{-1}, \ \ |\vec a|\geq 0
\end{split}
\end{equation}
\begin{equation}\notag
\begin{split}
&\lambda_{q+1}\|\nabla_{\vec a}\widetilde D_{t,q}^r\delta a_{I,q+1}\|_{0}+\|\nabla_{\vec a}\widetilde D_{t,q}^r \delta \theta_{I,q+1}\|_{0}+\|\nabla_{\vec a}\widetilde D_{t,q}^r \delta u_{I,q+1}\|_{0}\\
\lesssim&\ \lambda_{q+1}^{\frac12} \left(\frac{\lambda_{q+1}}{\lambda_q}\right)^{-1}\left(\frac{\lambda_{q+1}}{\lambda_q} \right)^{\frac{3}{2L}(|\vec a|+2-L)_+}\lambda_q^{|\vec a|}\delta_q^{\frac12}\tau_q^{-r}, \ \ |\vec a|\geq 0, \ \ r=0,1 \\
&\lambda_{q+1}\|\nabla_{\vec a}\widetilde D_{t,q}^2\delta a_{I,q+1}\|_{0}+\|\nabla_{\vec a}\widetilde D_{t,q}^2 \delta \theta_{I,q+1}\|_{0}+\|\nabla_{\vec a}\widetilde D_{t,q}^2 \delta u_{I,q+1}\|_{0}\\
\lesssim&\ \lambda_{q+1}^{\frac12} \left(\frac{\lambda_{q+1}}{\lambda_q}\right)^{-1}\left(\frac{\lambda_{q+1}}{\lambda_q} \right)^{\frac{3}{2L}(|\vec a|+2-L)_+}\lambda_q^{|\vec a|}\delta_q^{\frac12}\left(\lambda_{q-1}^{\frac32}\delta_{q-2}^{\frac12}\right)\tilde \tau_q^{-1}, \ \ |\vec a|\geq 0.
\end{split}
\end{equation}
\end{Lemma}

Thus we have
\begin{Lemma}\label{le-correct}
The estimates  
\begin{equation}\notag
\begin{split}
&\lambda_{q+1}\|\nabla_{\vec a}D_{t,q, \Pi,\mu}^r W_{I,q+1}\|_{0}+\|\nabla_{\vec a}D_{t,q, \Pi,\mu}^r \Lambda W_{I,q+1}\|_{0}+\|\nabla_{\vec a}D_{t,q, \Pi,\mu}^r \nabla^{\perp} W_{I,q+1}\|_{0}\\
\lesssim&\ \lambda_{q+1}^{\frac12+|\vec a|} \delta_q^{\frac12} \tau_q^{-r}, \ \ \forall \ \ |\vec a|\geq 0, \ \ r=0,1\\
&\lambda_{q+1}\|\nabla_{\vec a}D_{t,q, \Pi,\mu}^2 W_{I,q+1}\|_{0}+\|\nabla_{\vec a}D_{t,q, \Pi,\mu}^2 \Lambda W_{I,q+1}\|_{0}+\|\nabla_{\vec a}D_{t,q, \Pi,\mu}^2 \nabla^{\perp}W_{I,q+1}\|_{0}\\
\lesssim&\ \lambda_{q+1}^{\frac12+|\vec a|} \delta_q^{\frac12} \left(\lambda_q^{\frac32}\delta_{q-1}^{\frac12}\right)\tilde \tau_q^{-1}, \ \ \forall \ \ |\vec a|\geq 0,\\
\end{split}
\end{equation}
\begin{equation}\notag
\begin{split}
&\lambda_{q+1}\|\nabla_{\vec a}\widetilde D_{t,q}^r W_{I,q+1}\|_{0}+\|\nabla_{\vec a}\widetilde D_{t,q}^r \Lambda W_{I,q+1}\|_{0}+\|\nabla_{\vec a}\widetilde D_{t,q}^r \nabla^{\perp} W_{I,q+1}\|_{0}\\
\lesssim&\ \lambda_{q+1}^{\frac12+|\vec a|} \delta_q^{\frac12} \tau_{q-1}^{-r}, \ \ \forall \ \ |\vec a|\geq 0, \ \ r=0,1\\
&\lambda_{q+1}\|\nabla_{\vec a}\widetilde D_{t,q}^2 W_{I,q+1}\|_{0}+\|\nabla_{\vec a}\widetilde D_{t,q}^2 \Lambda W_{I,q+1}\|_{0}+\|\nabla_{\vec a}\widetilde D_{t,q}^2 \nabla^{\perp} W_{I,q+1}\|_{0}\\
\lesssim&\ \lambda_{q+1}^{\frac12+|\vec a|} \delta_q^{\frac12} \left(\lambda_{q-1}^{\frac32}\delta_{q-2}^{\frac12}\right)\tilde \tau_{q-1}^{-1}, \ \ \forall \ \ |\vec a|\geq 0
\end{split}
\end{equation}
hold for the correction terms.
\end{Lemma}

\medskip

\subsection{Estimates for second order anti-divergence operator}
\label{sec-anti}
We borrow the following lemmas from \cite{IM}.  
%{\color{red} include $\widetilde D_{t,q}=\partial_t+\widetilde u_q\cdot\nabla$ and $\widetilde D_{t,q+1}$ in the following two lemmas...}
\begin{Lemma}\label{le-aux-anti}
Let $Q:\mathbb R^2\to \mathbb R$ be Schwartz and satisfy
\[\max_{0\leq |\vec a|\leq |\vec b|\leq L}\lambda_{q+1}^{-|\vec b|+|\vec a|}\|h^{\vec a}\nabla^{\vec b} Q\|_{L^1(\mathbb R^2)}\leq \lambda^{-1}, \ \ \ |h|\sim \lambda_{q+1}^{-1}\]
for some $\lambda>0$. Then for smooth $U$ on $\mathbb T^2$ we have
\begin{equation}\notag
\begin{split}
\|\nabla_{\vec a}D_{t,q, \Pi,\mu}^r(Q*U)\|_{0}
\lesssim \lambda^{-1} \lambda_{q+1}^{|\vec a|}\left(\lambda_{q}^{\frac32}\delta_{q-1}^{\frac12}\right)^r \|U\|_{0},
\end{split}
\end{equation}
\begin{equation}\notag
\begin{split}
\|\nabla_{\vec a}D_{t,q+1, \Pi,\mu}^r(Q*U)\|_{0}
\lesssim \lambda^{-1} \lambda_{q+1}^{|\vec a|}\left(\lambda_{q+1}^{\frac32}\delta_{q}^{\frac12}\right)^r \|U\|_{0}.
\end{split}
\end{equation}
\end{Lemma}

\begin{Lemma}\label{le-anti-div}
Let $P_{\lambda_{q+1}}$ be a frequency projection operator localized to frequencies $\{\xi\in \widehat{\mathbb R}^2: \frac{1}{10}\lambda_{q+1} \leq |\xi|\leq 10\lambda_{q+1}\}$. Then for smooth $U$ on $\mathbb T^2$ we have
\begin{equation}\notag
\begin{split}
\|\nabla_{\vec a}D_{t,q, \Pi,\mu}^r(\mathcal R P_{\lambda_{q+1}}[U])\|_{C^0}
\lesssim  \lambda_{q+1}^{|\vec a|-2}\left(\lambda_{q}^{\frac32}\delta_{q-1}^{\frac12}\right)^r \|U\|_{C^0},
\end{split}
\end{equation}
\begin{equation}\notag
\begin{split}
\|\nabla_{\vec a}D_{t,q+1, \Pi,\mu}^r(\mathcal R P_{\lambda_{q+1}}[U])\|_{C^0}
\lesssim  \lambda_{q+1}^{|\vec a|-2}\left(\lambda_{q+1}^{\frac32}\delta_{q}^{\frac12}\right)^r \|U\|_{C^0}
\end{split}
\end{equation}
for $0\leq r\leq 1$ and $0\leq |\vec a|+r\leq L$.
\end{Lemma}

%\medskip

%\subsection{Basic estimates}
%\label{sec-basic}

\medskip

\subsection{The form of stress error}
\label{sec-stress-error}

It follows from the first equation of (\ref{g-q1}), noting $\widetilde u_q=\nabla^{\perp}\widetilde \eta_q$, $\widetilde \theta_q=\Lambda\widetilde \eta_q$ and (\ref{skip})
\begin{equation}\notag
\begin{split}
\nabla\cdot\nabla\cdot R_{q+1}=&\ \partial_t\Lambda W_{q+1}+\nabla\cdot(\widetilde u_q \Lambda W_{q+1})+\nabla\cdot(\nabla^{\perp} W_{q+1}\widetilde \theta_q)\\
&+\Lambda^{\gamma+1}W_{q+1}+\nabla\cdot\left(\nabla\cdot R_q-2\nabla^{\perp}W_{q+1}\Lambda W_{q+1}\right)\\
=&\ (\partial_t+ \widetilde u_{q-1}\cdot\nabla)\Lambda W_{q+1}+\nabla^{\perp} W_{q+1}\cdot\nabla\widetilde \theta_{q-1}+\Lambda^{\gamma+1}W_{q+1}\\
&+\nabla\cdot\left(\nabla\cdot \bar R_q-2\nabla^{\perp}W_{q+1}\Lambda W_{q+1}\right)+\nabla\cdot\nabla\cdot (R_q-\bar R_q)\\
=&: \nabla\cdot\nabla\cdot R_{T}+\nabla\cdot\nabla\cdot R_{N}+\nabla\cdot\nabla\cdot R_{D}+\nabla\cdot\nabla\cdot R_{O}+\nabla\cdot\nabla\cdot R_{M}
\end{split}
\end{equation}
with $R_T$, $R_N$, $R_D$, $R_O$ and $R_M$ representing the transport, Nash, dissipation, oscillation and mollification errors respectively. 

\medskip

\subsection{Transport error}
\label{sec-transport}

\begin{Lemma}\label{le-transport}
The estimates
\begin{equation}\notag
\begin{split}
%\|R_T\|_{C^0}\leq& \ C_0 \lambda_{q+1}^{-\frac32} \tau_{q}^{-1}\delta_q^{\frac12}\\
%\|D_{t,q}R_T\|_{C^0}\leq&\ C_0\lambda_{q+1}^{-\frac32} \tau_{q}^{-2}\delta_q^{\frac12}\\
\|R_{T}\|_{0}\lesssim& \ \lambda_{q+1}^{-\frac32} \tau_{q-1}^{-1}\delta_q^{\frac12}\\
\|D_{t,q,\Pi,\mu}R_T\|_{0}\lesssim&\ \lambda_{q+1}^{-\frac32}\tau_{q}^{-1}\tau_{q-1}^{-1}\delta_q^{\frac12}\\
\|\widetilde D_{t,q}R_T\|_{0}\lesssim&\ \lambda_{q+1}^{-\frac32} \tau_{q-1}^{-2}\delta_q^{\frac12}\\
%\|D_{t,q-1}R_{T}\|_{0}+\|\widetilde D_{t,q-1}R_T\|_{0}\lesssim&\ \lambda_{q+1}^{-\frac32} \tau_{q-1}^{-2}\delta_q^{\frac12}
\end{split}
\end{equation}
hold for the stress error $R_T$.
\end{Lemma}
\pf 
Note
\begin{equation}\notag
\nabla\cdot\nabla\cdot R_T=\left( \partial_t+\widetilde u_{q-1}\cdot\nabla \right)\Lambda W_{q+1}
%=&\ \Lambda\left( \left( \partial_t+u_q\cdot\nabla \right) W_{q+1}\right)+[u_q\cdot\nabla, \Lambda]W_{q+1}.
\end{equation}
where the right hand side has zero mean on $\mathbb T^2$. Hence there exists a symmetric trace-free tensor $R_T$ such that 
\[ R_T=\mathcal R\left( \partial_t+\widetilde u_{q-1}\cdot\nabla \right)\Lambda W_{q+1}\]
is well-defined. Noticing that $\left( \partial_t+\widetilde u_{q-1}\cdot\nabla \right)\Lambda W_{q+1}$ has frequency support near $\lambda_{q+1}$,
by Lemma \ref{le-correct} and Lemma \ref{le-anti-div}, we deduce
\begin{equation}\notag
\|R_{T}\|_{0}\lesssim \lambda_{q+1}^{-2} \lambda_{q+1}^{\frac12}\delta_q^{\frac12}\tau_{q-1}^{-1}\lesssim \lambda_{q+1}^{-\frac32} \tau_{q-1}^{-1}\delta_q^{\frac12}.
\end{equation}
Moreover 
\[ D_{t,q,\Pi,\mu}R_T=\mathcal RD_{t,q, \Pi, \mu}\left( \partial_t+\widetilde u_{q-1}\cdot\nabla \right)\Lambda W_{q+1}\]
\[ \widetilde D_{t,q}R_T=\mathcal R\widetilde D_{t,q}\left( \partial_t+\widetilde u_{q-1}\cdot\nabla \right)\Lambda W_{q+1}\]
and hence
\begin{equation}\notag
\begin{split}
\|D_{t,q, \Pi, \mu}R_T\|_{0}\lesssim&\ \lambda_{q+1}^{-2} \lambda_{q+1}^{\frac12}\delta_q^{\frac12}\tau_{q}^{-1}\tau_{q-1}^{-1}\lesssim \lambda_{q+1}^{-\frac32} \tau_{q}^{-1}\tau_{q-1}^{-1}\delta_q^{\frac12},\\
\|\widetilde D_{t,q}R_T\|_{0}\lesssim&\ \lambda_{q+1}^{-2} \lambda_{q+1}^{\frac12}\delta_q^{\frac12}\tau_{q-1}^{-1}\tau_{q-1}^{-1}\lesssim \lambda_{q+1}^{-\frac32}\tau_{q-1}^{-2}\delta_q^{\frac12}.
\end{split}
\end{equation}
%The estimates for $\|\widetilde D_{t,q}R_T\|_{0}$ and $\|\widetilde D_{t,q-1}R_T\|_{0}$ can be obtained similarly.
\cbdu

\medskip

\subsection{Nash error}
\label{sec-Nash}

\begin{Lemma}\label{le-Nash}
The Nash error $R_N$ satisfies the following estimates
\begin{equation}\notag
\begin{split}
\|R_{N}\|_{0}\lesssim&\ \lambda_{q+1}^{-\frac32} \lambda_{q-1}^{\frac32}\delta_q^{\frac12}\delta_{q-2}^{\frac12} \\
\|D_{t,q,\Pi,\mu}R_N\|_{0}\lesssim&\ \tau_{q}^{-1} \lambda_{q+1}^{-\frac32} \lambda_{q-1}^{\frac32}\delta_q^{\frac12}\delta_{q-2}^{\frac12}\\
\|\widetilde D_{t,q}R_N\|_{0}\lesssim&\ \tau_{q-1}^{-1} \lambda_{q+1}^{-\frac32} \lambda_{q-1}^{\frac32}\delta_q^{\frac12}\delta_{q-2}^{\frac12}.
%\|D_{t,q-1}R_{N}\|_{0}\lesssim&\ \tau_{q-1}^{-1} \lambda_{q+1}^{-\frac32} \lambda_{q-1}^{\frac32}\delta_q^{\frac12}\delta_{q-2}^{\frac12}.
\end{split}
\end{equation}
\end{Lemma}
\pf
Since $\nabla\cdot\nabla\cdot R_N=\nabla^{\perp}W_{q+1}\cdot\nabla\widetilde \theta_{q-1}$, the right hand side of which also has zero mean on $\mathbb T^2$, there is a symmetric trace-free tensor $R_N$ defined by
\begin{equation}\notag
R_N=\mathcal R \left(\nabla^{\perp}W_{q+1}\cdot\nabla\widetilde \theta_{q-1}\right).
\end{equation}
We also note $\nabla^{\perp}W_{q+1}\cdot\nabla\widetilde \theta_{q-1}$ has frequency support near $\lambda_{q+1}$.
In view of the induction estimate (\ref{induct-theta}), Lemma \ref{le-correct} and Lemma \ref{le-anti-div} we have
\begin{equation}\notag
\begin{split}
\|R_{N}\|_{0}\lesssim &\ \lambda_{q+1}^{-2} \|\nabla \widetilde \theta_{q-1}\|_{0}\|\nabla^{\perp}W_{q+1}\|_{0}\\
\lesssim &\ \lambda_{q+1}^{-2} \sum_{I\in \mathbb Z\times F}\lambda_{q-1}\|\widetilde\theta_{q-1}\|_{0}\lambda_{q+1}\|a_{I,q+1}\|_{0}\\
\lesssim&\ \lambda_{q+1}^{-1} \lambda_{q-1}^{\frac32}\delta_{q-2}^{\frac12}\delta_q^{\frac12}\lambda_{q+1}^{-\frac12}.
\end{split}
\end{equation}
To estimate $\widetilde D_{t,q} R_N$, we have
\begin{equation}\notag
\widetilde D_{t,q} \left(\nabla^{\perp}W_{q+1}\cdot\nabla \widetilde\theta_{q-1} \right)=\widetilde D_{t,q} \left(\nabla^{\perp}W_{q+1}\right)\cdot\nabla \widetilde\theta_{q-1}+ \nabla^{\perp}W_{q+1}\cdot \widetilde D_{t,q}\nabla \widetilde\theta_{q-1}.
\end{equation}
We observe
\begin{equation}\notag
\begin{split}
\widetilde D_{t,q}\nabla \widetilde\theta_{q-1}=&\ (\partial_t+\widetilde u_{q-1}\cdot\nabla)\nabla \widetilde\theta_{q-1}\\
=&\ (\partial_t+\nabla^{\perp}\Pi_{q-1}\cdot\nabla)\nabla \Lambda\Pi_{q-1}-(\partial_t+\nabla^{\perp}\Pi_{q-1}\cdot\nabla)\nabla \Lambda\mu_{q-1}\\
&-(\partial_t+\nabla^{\perp}\mu_{q-1}\cdot\nabla)\nabla \Lambda\Pi_{q-1}+(\partial_t+\nabla^{\perp}\mu_{q-1}\cdot\nabla)\nabla \Lambda\mu_{q-1}
\end{split}
\end{equation}
and hence by the inductive estimates (\ref{induct-mu-t})
\begin{equation}\notag
\|\widetilde D_{t,q}\nabla \widetilde\theta_{q-1}\|_{0}\lesssim \lambda_{q-1}^3\delta_{q-2}.
\end{equation}
By Lemma \ref{le-correct} and Lemma \ref{le-amplitude-main},
\begin{equation}\notag
\begin{split}
\|\widetilde D_{t,q} R_N\|_{0}\lesssim&\ \lambda_{q+1}^{-2}\|\widetilde D_{t,q} \nabla^{\perp}W_{q+1}\|_{0}\|\nabla\widetilde\theta_{q-1}\|_{0}
+\lambda_{q+1}^{-2}\|\nabla^{\perp}W_{q+1}\|_{0}\|\widetilde D_{t,q}\nabla\widetilde\theta_{q-1}\|_{0}\\
\lesssim&\ \lambda_{q+1}^{-2}\lambda_{q+1}^{\frac12}\delta_q^{\frac12}\tau_{q-1}^{-1}\lambda_{q-1}^{\frac32}\delta_{q-2}^{\frac12}+\lambda_{q+1}^{-2}\lambda_{q+1}^{\frac12}\delta_q^{\frac12}\lambda_{q-1}^{3}\delta_{q-2}\\
\lesssim&\ \lambda_{q+1}^{-\frac32}\lambda_{q-1}^{\frac32}\delta_q^{\frac12}\delta_{q-2}^{\frac12}\tau_{q-1}^{-1}
\end{split}
\end{equation}
where we used the fact $\lambda_{q-1}^{\frac32}\delta_{q-2}^{\frac12}\lesssim \tau_{q-1}^{-1}$.

The estimate of $\|D_{t,q,\Pi,\mu}R_N\|_{0}$ can be obtained similarly by noticing that the cost of $D_{t,q,\Pi}$ and $D_{t,q,\mu}$ is $\tau_q^{-1}$.

\cbdu

\medskip

\subsection{Dissipation error}
\label{sec-dissipation}

\begin{Lemma}\label{le-dis}
The dissipation error $R_D$ has the estimates
\begin{equation}\notag
\begin{split}
\|R_{D}\|_{0}\lesssim&\ \lambda_{q+1}^{\gamma-\frac32} \delta_q^{\frac12}\\
\|D_{t,q,\Pi,\mu}R_{D}\|_{0}\lesssim&\ \lambda_{q+1}^{\gamma-\frac32} \delta_q^{\frac12}\tau_{q}^{-1}\\
\|\widetilde D_{t,q}R_{D}\|_{0}\lesssim&\ \lambda_{q+1}^{\gamma-\frac32} \delta_q^{\frac12}\tau_{q-1}^{-1}.
\end{split}
\end{equation}
\end{Lemma}
\pf
It follows from 
\[\nabla\cdot\nabla\cdot R_D=\Lambda^{\gamma+1}W_{q+1}\]
that $R_D=\mathcal R\Lambda^{\gamma+1}W_{q+1}$. Applying Lemma \ref{le-correct} and Lemma \ref{le-anti-div} yields
\begin{equation}\notag
\|R_{D}\|_{0}\lesssim \lambda_{q+1}^{-2} \lambda_{q+1}^{\gamma+1}\|W_{q+1}\|_{0}
\lesssim \lambda_{q+1}^{\gamma-\frac32} \delta_q^{\frac12}.
\end{equation}
While 
\[\widetilde D_{t,q}R_D=\mathcal R \Lambda^{\gamma+1}\widetilde D_{t,q}W_{q+1}+[\widetilde D_{t,q}, \mathcal R\Lambda^{\gamma+1}]W_{q+1},\]
we deduce from Lemma \ref{le-correct}, Lemma \ref{le-anti-div}, the commutator estimate from \cite{BSV} and the inductive estimate (\ref{induct-u1})
\begin{equation}\notag
\begin{split}
\|\widetilde D_{t,q}R_D\|_0\lesssim&\ \lambda_{q+1}^{-2+\gamma+1}\lambda_{q+1}^{-\frac12}\delta_q^{\frac12}\tau_{q-1}^{-1}
+\lambda_{q+1}^{-2+\gamma+1}\lambda_{q+1}^{-\frac12}\delta_q^{\frac12}\|\nabla\widetilde u_{q-1}\|_0\\
\lesssim&\ \lambda_{q+1}^{-2+\gamma+1}\lambda_{q+1}^{-\frac12}\delta_q^{\frac12}\tau_{q-1}^{-1}
+\lambda_{q+1}^{-2+\gamma+1}\lambda_{q+1}^{-\frac12}\delta_q^{\frac12}\lambda_{q-1}^{\frac32}\delta_{q-1}^{\frac12}\\
\lesssim&\ \lambda_{q+1}^{\gamma-\frac32}\delta_q^{\frac12}\tau_{q-1}^{-1}.
\end{split}
\end{equation}
The terms $D_{t,q,\Pi}R_{D}$ and $D_{t,q,\mu}R_{D}$ can be handled analogously. 

\cbdu

\medskip

\subsection{Oscillation error}
\label{sec-osc}
\begin{Lemma}\label{le-osc}
The estimates
\begin{equation}\notag
\begin{split}
\|R_{O}\|_{0}\lesssim&\ \tau_{q-1}\lambda_{q-1}^{\frac32}\delta_q\delta_{q-2}^{\frac12} \\
\|D_{t,q,\Pi, \mu}R_{O}\|_{0}\lesssim&\ \delta_q\tau_{q}^{-1}\\
\|\widetilde D_{t,q}R_{O}\|_{0}\lesssim&\ \delta_q\tau_{q-1}^{-1}
\end{split}
\end{equation}
are valid for the oscillation error $R_O$.
\end{Lemma}
\pf
Note $\nabla\cdot\nabla\cdot R_O$ can be decomposed as 
\begin{equation}\notag
\begin{split}
\nabla\cdot\nabla\cdot R_O=&\left(\nabla\cdot\nabla\cdot\bar R_q+2\nabla\cdot \sum_{I=(k,v)\in \mathbb Z\times F}\nabla^{\perp}W_{I,q+1} \Lambda W_{\bar I, q+1}\right)\\
&+ \sum_{I,J\in \mathbb Z\times F, J\neq \bar I}\left(\nabla^{\perp} W_{I,q+1}\nabla\Lambda W_{J, q+1} +\nabla^{\perp} W_{J, q+1}\nabla\Lambda W_{I, q+1} \right)\\
=&: \nabla\cdot\nabla\cdot R_S+\nabla\cdot\nabla\cdot R_H.
\end{split}
\end{equation}
Through the proof of Lemma \ref{le-decomp} we know
\[R_{S}= \sum_{I=(k,v)\atop \in \mathbb Z\times F}\delta B_I.\]
On the other hand we have
\begin{equation}\notag
R_{H}= \mathcal R \sum_{I,J\in \mathbb Z\times F, J\neq \bar I}\left(\nabla^{\perp} W_{I,q+1}\nabla\Lambda W_{J, q+1} +\nabla^{\perp} W_{J, q+1}\nabla\Lambda W_{I, q+1} \right).
\end{equation}
Recall
\[W_{I,q+1}= P_{I,\lambda_{q+1}}\left(a_{I,q+1}(x,t)e^{i\lambda_{q+1}\xi_I(x,t)} \right)=e^{i\lambda_{q+1}\xi_I(x,t)}(a_{I,q+1}+\delta a_{I,q+1})\]
and 
\begin{equation}\notag
\begin{split}
\Lambda W_{I,q+1}=&\ e^{i\lambda_{q+1}\xi_I(x,t)}(\theta_{I,q+1}+\delta \theta_{I,q+1})\\
\nabla\Lambda W_{I,q+1}=&\ i\lambda_{q+1}e^{i\lambda_{q+1}\xi_I(x,t)}\nabla\xi_I(\theta_{I,q+1}+\delta \theta_{I,q+1})\\
&+e^{i\lambda_{q+1}\xi_I(x,t)}\nabla(\theta_{I,q+1}+\delta \theta_{I,q+1})\\
\nabla^{\perp} W_{I,q+1} =&\ i\lambda_{q+1} e^{i\lambda_{q+1}\xi_I(x,t)}\nabla^{\perp} \xi_I (a_{I,q+1}+\delta a_{I,q+1})\\
&+e^{i\lambda_{q+1}\xi_I(x,t)}\nabla^{\perp}(a_{I,q+1}+\delta a_{I,q+1})
\end{split}
\end{equation}
with
\[\theta_{I,q+1}=\lambda_{q+1}|\nabla \xi_I|a_{I,q+1}.\]
Therefore we have
\begin{equation}\notag
\begin{split}
&\nabla^{\perp} W_{I,q+1}\nabla\Lambda W_{J,q+1}+\nabla^{\perp} W_{J,q+1}\nabla\Lambda W_{I,q+1}\\
=&-\lambda_{q+1}^2e^{i\lambda_{q+1}(\xi_I+\xi_J)}\left(\nabla^{\perp}\xi_I\nabla\xi_J a_{I,q+1}\theta_{J, q+1}+\nabla^{\perp}\xi_J\nabla\xi_I a_{J,q+1}\theta_{I, q+1}\right)\\
&-\lambda_{q+1}^2e^{i\lambda_{q+1}(\xi_I+\xi_J)}\left(\nabla^{\perp}\xi_I\nabla\xi_J \delta a_{I,q+1}\theta_{J, q+1}+\nabla^{\perp}\xi_J\nabla\xi_I \delta a_{J,q+1}\theta_{I, q+1}\right)\\
&-\lambda_{q+1}^2e^{i\lambda_{q+1}(\xi_I+\xi_J)}\left(\nabla^{\perp}\xi_I\nabla\xi_J \delta \theta_{I,q+1}(a_{J, q+1}+\delta a_{J, q+1})\right.\\
&\left. \ \ \ \ \ \ \ \ \ \ \ \ \ \ \ \ \ \ \ \ \ \ \ \ \ \ \ \ +\nabla^{\perp}\xi_J\nabla\xi_I \delta \theta_{J,q+1}(a_{I, q+1}+\delta a_{I, q+1})\right)\\
&+i\lambda_{q+1}e^{i\lambda_{q+1}(\xi_I+\xi_J)} \left(\nabla\xi_J(a_{I,q+1}+\delta a_{I,q+1})(\theta_{J,q+1}+\delta \theta_{J,q+1}) \right.\\
&\left. \ \ \ \ \ \ \ \ \ \ \ \ \ \ \ \ \ \ \ \ \ \ \ \ \ \ \ \ +\nabla\xi_I(a_{J,q+1}+\delta a_{J,q+1})(\theta_{I,q+1}+\delta \theta_{I,q+1})\right)\\
&+e^{i\lambda_{q+1}(\xi_I+\xi_J)} \left(\nabla^{\perp}(a_{I,q+1}+\delta a_{I,q+1})\nabla(\theta_{J,q+1}+\delta \theta_{J,q+1}) \right.\\
&\left. \ \ \ \ \ \ \ \ \ \ \ \ \ \ \ \ \ \ \ \ \ \ +\nabla^{\perp}(a_{J,q+1}+\delta a_{J,q+1})\nabla(\theta_{I,q+1}+\delta \theta_{I,q+1})\right)\\
=&: K_1+K_2+K_3+K_4+K_5
\end{split}
\end{equation}
where $K_1$ is the leading order term. We note
\begin{equation}\notag
\begin{split}
&\nabla^{\perp}\xi_I\nabla\xi_J a_{I,q+1}\theta_{J, q+1}+\nabla^{\perp}\xi_J\nabla\xi_I a_{J,q+1}\theta_{I, q+1}\\
=&\lambda_{q+1}a_{I,q+1}a_{J, q+1}\left(\nabla^{\perp}\xi_I\nabla\xi_J |\nabla\xi_J| +\nabla^{\perp}\xi_J\nabla\xi_I |\nabla\xi_I| \right)\\
=&\lambda_{q+1}a_{I,q+1}a_{J, q+1}\left(\nabla^{\perp}\xi_I\nabla\xi_J (|\nabla\xi_J|-|\nabla\xi_{J,in}|) +\nabla^{\perp}\xi_J\nabla\xi_I (|\nabla\xi_I| -|\nabla\xi_{I,in}|) \right)\\
\end{split}
\end{equation}
by using the fact that $\nabla\xi_{I,in}=\nabla\xi_{J,in}$ and 
\[\nabla^{\perp}\xi_I\nabla\xi_J+\nabla^{\perp}\xi_J\nabla\xi_I=0.\]
Therefore we have
\begin{equation}\notag
\begin{split}
\|\mathcal RK_1\|_{0}\lesssim&\ \lambda_{q+1}^{-2+3} \sum \|a_{I,q+1}\|_{0}^2\left||\nabla\xi_I|-|\nabla\xi_{I,in}|\right|\\
\lesssim&\ \lambda_{q+1}^{-2+3}(\lambda_{q+1}^{-\frac12}\delta_q^{\frac12})^2\lambda_{q-1}\tau_{q-1}\|\widetilde u_{q-1}\|_{0}\\
\lesssim&\ \lambda_{q+1}^{-2+3}(\lambda_{q+1}^{-\frac12}\delta_q^{\frac12})^2\lambda_{q-1}\tau_{q-1}\lambda_{q-1}^{\frac12}\delta_{q-2}^{\frac12}\\
\lesssim&\ \tau_{q-1}\lambda_{q-1}^{\frac32}\delta_q\delta_{q-2}^{\frac12}.
\end{split}
\end{equation}
It is clear that $K_j$'s are lower order terms for $2\leq j\leq 5$. Hence the estimate details for these terms are omitted and we claim
\[\|R_H\|_{0}\lesssim \tau_{q-1}\lambda_{q-1}^{\frac32}\delta_q\delta_{q-2}^{\frac12}.\]

We further estimate $\widetilde D_{t,q}K_1$. Straightforward computation shows
\begin{equation}\notag
\begin{split}
&\widetilde D_{t,q}\left(\nabla^{\perp}\xi_I\nabla\xi_J a_{I,q+1}\theta_{J, q+1}+\nabla^{\perp}\xi_J\nabla\xi_I a_{J,q+1}\theta_{I, q+1}\right)\\
=& \left[(\widetilde D_{t,q}\nabla^{\perp}\xi_I)\nabla\xi_J a_{I,q+1}\theta_{J, q+1}+(\widetilde D_{t,q}\nabla^{\perp}\xi_J)\nabla\xi_I a_{J,q+1}\theta_{I, q+1}\right]\\
&+\left[(\nabla^{\perp}\xi_I(\widetilde D_{t,q}\nabla\xi_J) a_{I,q+1}\theta_{J, q+1}+\nabla^{\perp}\xi_J(\widetilde D_{t,q}\nabla\xi_I) a_{J,q+1}\theta_{I, q+1}\right]\\
&+\left[(\nabla^{\perp}\xi_I\nabla\xi_J (\widetilde D_{t,q}a_{I,q+1})\theta_{J, q+1}+\nabla^{\perp}\xi_J\nabla\xi_I (\widetilde D_{t,q}a_{J,q+1})\theta_{I, q+1}\right]\\
&+\left[(\nabla^{\perp}\xi_I\nabla\xi_J a_{I,q+1}(\widetilde D_{t,q}\theta_{J, q+1})+\nabla^{\perp}\xi_J\nabla\xi_I a_{J,q+1}(\widetilde D_{t,q}\theta_{I, q+1})\right].
%=&\ \widetilde K_1+\widetilde K_2+\widetilde K_3+\widetilde K_4.
\end{split}
\end{equation}
Since $\theta_{I,q+1}=\lambda_{q+1}|\nabla\xi_I|a_{I,q+1}$, by Lemma \ref{le-phase} and Lemma \ref{le-amplitude-main}, we have
\begin{equation}\notag
\begin{split}
\|\theta_{I,q+1}\|_{0}\lesssim&\ \lambda_{q+1}^{\frac12}\delta_q^{\frac12},\\
\|\widetilde D_{t,q}\theta_{I,q+1}\|_{0}\lesssim&\ \lambda_{q+1} \|\widetilde D_{t,q}\nabla\xi_I\|_{0}\|a_{I,q+1}\|_{0}+\lambda_{q+1}\|\widetilde D_{t,q}a_{I,q+1}\|_{0}\\
\lesssim&\ \lambda_{q+1}^{\frac12}\delta_q^{\frac12}\lambda_{q-1}^{\frac32}\delta_{q-2}^{\frac12}+\lambda_{q+1}^{\frac12}\delta_q^{\frac12}\tau_{q-1}^{-1}\\
\lesssim&\ \lambda_{q+1}^{\frac12}\delta_q^{\frac12}\tau_{q-1}^{-1}
\end{split}
\end{equation}
where we used the fact $\lambda_{q-1}^{\frac32}\delta_{q-2}^{\frac12}\lesssim \tau_{q-1}^{-1}$.
Hence combining with Lemma \ref{le-phase} and Lemma \ref{le-amplitude-main}, it follows 
\begin{equation}\notag
\begin{split}
&\|\mathcal R \widetilde D_{t,q}\left(\nabla^{\perp}\xi_I\nabla\xi_J a_{I,q+1}\theta_{J, q+1}+\nabla^{\perp}\xi_J\nabla\xi_I a_{J,q+1}\theta_{I, q+1}\right)\|_{0}\\
\lesssim &\ \lambda_{q+1}^{-2+2} \lambda_{q-1}^{\frac32}\delta_{q-2}^{\frac12}\lambda_{q+1}^{-\frac12}\delta_q^{\frac12}\lambda_{q+1}^{\frac12}\delta_q^{\frac12}+\lambda_{q+1}^{-2+2}\lambda_{q+1}^{-\frac12}\delta_q^{\frac12}\tau_{q-1}^{-1}\lambda_{q+1}^{\frac12}\delta_q^{\frac12}\\
\lesssim&\ \delta_q\tau_{q-1}^{-1}.
\end{split}
\end{equation}
Other terms in $\widetilde D_{t,q} R_H$ can be analyzed similarly and hence 
\[\|\widetilde D_{t,q} R_H\|_0 \lesssim\delta_q\tau_{q-1}^{-1}.\]
On the other hand, we have from Lemma \ref{le-error1-micro} that
\begin{equation}\notag
\begin{split}
\|R_{S}\|_{0}\lesssim&\ \lambda_{q+1}^{-2}(\lambda_{q+1}^{-1}\lambda_q)\delta_q\lesssim \tau_{q-1}\lambda_{q-1}^{\frac32}\delta_q\delta_{q-2}^{\frac12},\\
\|\widetilde D_{t,q}R_{S}\|_{0}\lesssim&\ \lambda_{q+1}^{-2}(\lambda_{q+1}^{-1}\lambda_q)\delta_q \tau_{q-1}^{-1}\lesssim \delta_q\tau_{q-1}^{-1}.
\end{split}
\end{equation}
Again, the estimate for $\|D_{t,q,\Pi,\mu} R_H\|_0$ can be obtained analogously. It concludes the proof.

\cbdu

\medskip

\subsection{Mollification error}
\label{sec-molli-error}
\begin{Lemma}\label{le-molli-error}
The mollification error $R_M$ satisfies the estimates
\begin{equation}\notag
\begin{split}
\|R_{M}\|_{0}\lesssim&\ \left(\frac{\lambda_{q}}{\lambda_{q-1}}\right)^{-\frac32}\delta_q \\
\|D_{t,q,\Pi, \mu}R_{M}\|_{0}\lesssim&\ \lambda_{q}^{\frac32}\delta_{q-1}^{\frac12}\delta_q\\
\|\widetilde D_{t,q}R_{M}\|_{0}\lesssim&\ \lambda_{q-1}^{\frac32}\delta_{q-2}^{\frac12}\delta_q.
\end{split}
\end{equation}
\end{Lemma}
\pf
As discussed in Section \ref{sec-molli}, for mollified stress $\bar R_q$ via the flow map of $\partial_t+u_q\cdot\nabla$, estimates in Lemma \ref{le-molli} hold. Thus for mollified stress $\bar R_q$ via the flow map of $\partial_t+\widetilde u_q\cdot\nabla=\partial_t+\widetilde u_{q-1}\cdot\nabla$, we have
\begin{equation}\notag
\|R_M\|_0\lesssim \left(\left( \lambda_{q-1}^{\frac32}\delta_{q-2}^{\frac12}\right)\tilde \tau_{q-1}+l_{q-1}^L\lambda_{q-1}^L \right) \delta_q \lesssim \left(\frac{\lambda_{q}}{\lambda_{q-1}}\right)^{-\frac32}\delta_q
\end{equation}
where we used the choice of mollification scales (\ref{scale-molli-xt}) with $q$ replaced by $q-1$. For the material derivative, we only show the estimate of $\widetilde D_{t,q}R_M$ as the other one can be obtained similarly. We infer from the inductive estimate (\ref{induct-mat1}) and Lemma \ref{le-molli} for the flow map associated with $\widetilde u_q=\widetilde u_{q+1}$
\begin{equation}\notag
\|\widetilde D_{t,q}R_M\|_0\lesssim \|\widetilde D_{t,q}R_q\|_0+\|\widetilde D_{t,q}\bar R_q\|_0\lesssim \lambda_{q-1}^{\frac32}\delta_{q-2}^{\frac12}\delta_q.
\end{equation}
\cbdu

\subsection{Concluding the proof of Proposition \ref{prop-main}}
We need to show that the solution $(\Pi_{q+1}, \mu_{q+1}, \widetilde R_{q+1}, R_{q+1})$ satisfies the estimates (\ref{induct-theta})-(\ref{induct-freq}) with $q$ replaced by $q+1$. Since $\mu_{q+1}=\mu_q+W_{q+1}$ and $\Pi_{q+1}=\Pi_q-W_{q+1}$, it follows from (\ref{induct-theta}) and Lemma \ref{le-correct}
\begin{equation}\notag
\begin{split}
\|\mu_{q+1}\|_{N+1}+\|\Pi_{q+1}\|_{N+1}\leq&\ \|\mu_{q}\|_{N+1}+\|\Pi_{q}\|_{N+1}+2\|W_{q+1}\|_{N+1}\\
\leq&\ M\lambda_q^{N+\frac12}\delta_{q-1}^{\frac12}
+\lambda_{q+1}^{N+\frac12}\delta_{q}^{\frac12}\\
\leq&\ M\lambda_{q+1}^{N+\frac12}\delta_{q}^{\frac12}.
\end{split}
\end{equation}
Thus (\ref{induct-theta}) is true for $q$ replaced by $q+1$. 

Combining Lemmas \ref{le-transport}-\ref{le-molli-error} we have
\begin{equation}\label{Rq1-final}
\begin{split}
\|R_{q+1}\|_0\lesssim&\ \lambda_{q+1}^{-\frac32}\tau_{q-1}^{-1}\delta_q^{\frac12}+\lambda_{q+1}^{-\frac32}\lambda_{q-1}^{\frac32}\delta_q^{\frac12}\delta_{q-2}^{\frac12}+\tau_{q-1}\lambda_{q-1}^{\frac32}\delta_q\delta_{q-2}^{\frac12}\\
&+\lambda_{q+1}^{\gamma-\frac32}\delta_q^{\frac12}+\left(\frac{\lambda_{q}}{\lambda_{q-1}}\right)^{-\frac32}\delta_q.
\end{split}
\end{equation}
As discussed in Section \ref{sec-heuristics}, the first three terms in (\ref{Rq1-final}) are $\lesssim \delta_{q+1}$ provided $\beta<1$.  While 
\[\lambda_{q+1}^{\gamma-\frac32}\delta_q^{\frac12}\lesssim \delta_{q+1} \ \ \ \ \mbox{for} \ \ \beta<\frac{b-2b(\gamma-1)}{2b-1},\]
\[\left(\frac{\lambda_{q}}{\lambda_{q-1}}\right)^{-\frac32}\delta_q\lesssim \delta_{q+1} \ \ \ \ \mbox{for} \ \ \beta<\frac32.\]
Combining with (\ref{alpha}) we claim that $\|R_{q+1}\|_0\leq \delta_{q+1}$ and hence $\|R_{q+1}\|_N\leq \lambda_{q+1}^N\delta_{q+1}$ since $\widehat R_{q+1}$ is supported near $\lambda_{q+1}$, that is (\ref{induct-R}) is satisfied with $q$ replaced by $q+1$, providing $-\frac12\leq \alpha<0$ and $0\leq \gamma<1-\alpha$.

Since 
\begin{equation}\notag
\begin{split}
(\partial_t+\nabla^{\perp}\Pi_{q+1}\cdot\nabla)\Lambda\mu_{q+1}=&\ (\partial_t+\nabla^{\perp}\Pi_{q}\cdot\nabla)\Lambda\mu_q-\nabla^{\perp}W_{q+1}\cdot\nabla \Lambda\mu_q\\
&-\nabla^{\perp}W_{q+1}\cdot\nabla \Lambda W_{q+1}+\nabla^{\perp}\Pi_{q}\cdot\nabla \Lambda W_{q+1},
\end{split}
\end{equation}
applying (\ref{induct-theta}), (\ref{induct-mu-t}) and Lemma \ref{le-correct} gives
\[\|(\partial_t+\nabla^{\perp}\Pi_{q+1}\cdot\nabla)\Lambda\mu_{q+1}\|_N\lesssim \lambda_{q+1}^{N+2}\delta_q. \]
Other terms in (\ref{induct-mu-t}) with $q$ replaced by $q+1$ can be verified similarly. 

Estimate (\ref{induct-R-t}) with $q$ replaced by $q+1$ follows from Lemmas \ref{le-transport} -\ref{le-molli-error}. The frequency support condition (\ref{induct-freq}) and time support property (\ref{time-supp}) with $q$ replaced by $q+1$ is certainly true in view of the construction of $W_{q+1}$. 

Moreover, we note
\begin{equation}\notag
\partial_t \Lambda^{\frac12}W_{q+1}=\Lambda^{\frac12}D_{t,q} W_{q+1}-\Lambda^{\frac12}(u_q\cdot\nabla W_{q+1})
\end{equation}
which implies 
\[\|\Lambda^{\frac12}W_{q+1}\|_{C_t^1 C_x^0}\lesssim \lambda_{q+1}^{\frac12}\lambda_q^{\frac12}\delta_{q-1}^{\frac12}\lambda_{q+1}^{\frac12}\delta_q^{\frac12}\lesssim \lambda_{q+1}\lambda_q^{\frac12}\delta_{q-1}^{\frac12}\delta_q^{\frac12},\]
and hence by interpolation
\[\|\Lambda^{\frac12}W_{q+1}\|_{C_t^\zeta C_x^0}\lesssim (\lambda_{q+1}\lambda_q^{\frac12}\delta_{q-1}^{\frac12})^{\zeta}\delta_q^{\frac12}\sim \lambda_q^{(b+\frac12-\frac{1}{2b}\beta)\zeta-\frac{1}{2}\beta}. \]
Therefore, to have $C^\zeta$ regularity in time, it is necessary to impose
\[(b+\frac12-\frac{1}{2b}\beta)\zeta-\frac{1}{2}\beta<0\]
which indicates
\[\zeta<\frac{\beta}{2b+1-\frac{1}{b}\beta}<\frac12\]
when choosing $b=1^+$ and $\beta<1$.

\bigskip

\section{Proof of Theorem \ref{thm}}
\label{sec-proof}
With the flexibility of allowing external forcing in both equations of (\ref{pm-q}), we can find an initial tuple $(\Pi_0,\mu_0,\widetilde R_0, R_0)$ satisfying the system (\ref{pm-q}) and (\ref{induct-theta})-(\ref{induct-freq}) at level 0, such that $\mu_0\not\equiv 0$. Applying the inductive Proposition \ref{prop-main} generates a sequence $\{(\Pi_q,\mu_q,\widetilde R_q, R_q)\}$ satisfying  (\ref{pm-q}) and (\ref{induct-theta})-(\ref{induct-freq}). With the estimates (\ref{induct-theta}) and (\ref{induct-R}), taking the limit $q\to \infty$ yields a limit solution $(\Pi, \mu, \widetilde R, 0)$ of (\ref{pm3}) with $\mu\not\equiv 0$. Equivalently, there are at least two solutions to (\ref{sqg}) with external forcing $f=\nabla\cdot\nabla\cdot \widetilde R$. In the end we note from (\ref{g-q1})- (\ref{g-q2}) that $f\in C_t^0C_x^{2\alpha-1}$.

%\bigskip

%\section*{Acknowledgement}
%The authors would like to express their gratitude to ... 

%\Endrefs
\end{document}